# Avoiding Intersections of Dragon Curves

REIMUND ALBERS, ZONGYI GUO, HUAIYI GUO[1]

Abstract: This article proves that there are no self-intersections in the dragon curve when the unfolding angle is greater than 98.195°. This is shown by constructing a hull for the dragon curve that is mapped onto itself by the generating mappings for the dragon curve. The treatment is purely geometric. The proof is supplemented by a conjecture that reduces the boundary for the unfolding angle to 96.241°
.

**1. Introduction**

The paperfolding polygon, when unfolded by exactly 90°, exhibits no intersections, meaning that no segment of the broken line doubles back on itself [1]. The question if intersections are avoided or not when the unfolding angle exceeds 90° first came to attention [2] in 1990, stemming from M. Mendès France's intuitive assumption, that no intersections occur with unfolding angles greater than 90°. Computer experiments conducted around 1995 by H.-O. Peitgen and S. Michaelis provided indications of intersections occurring for unfolding angles between 90° and 95°. In his 2011 book [1], D. Knuth noted that he "*noticed in 1969 that 95°-angle folds would lead to paths that cross themselves*". In 2006, Albers presented a proof [3] demonstrating the existence of crossing segments in iteration 10 for angles ranging between 90° and 95.126°. It is additionally demonstrated [3] that certain vertices of the paper-folding curve coincide when the unfolding angle is approximately 93.06°. In 2014, S. Tabachnikov raises in "Dragon curves revisited" [4] the question: "*…, what is the value of the critical angle for which the curve starts to cross itself, and where does this self-crossing occur?*" In 2018, J.-P. Allouche, M. Mendès France, and G. Skordev provide a comprehensive overview [5] of the problem and its solutions to date. In 2022, Y. Kamiya constructs a sequence of angles [6] with the limit of 93.7912°, where two vertices of the paper-folding curve precisely coincide. Finally, in January 2024, S. Akiyama, Y. Kamiya, and F. Wen establish that no intersections occur [7] when the unfolding angle exceeds 99.3438°. Their work and the present article address the same problem and use the same big picture. While their approach employs complex calculations focused on the limit Dragon curve, this article presents a purely geometric analysis of the finite paperfolding polygons.

The current findings are summarized in Figure 1.

[1]*Universität Bremen. MZH Bibliothekstr. 5, 28359 Bremen, Germany*
*Gymnasium Horn, Vorkampsweg 97, 28359 Bremen, Germany*
*E-mail addresses: ralbers@uni-bremen.de, zy.guo@dxfund.org, hy. guo@dxfund.org*

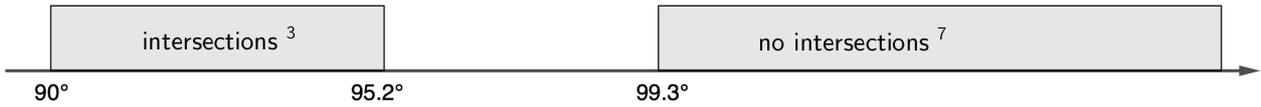

Figure 1: The present knowledge

In this article the gap of uncertainty is narrowed as follows.

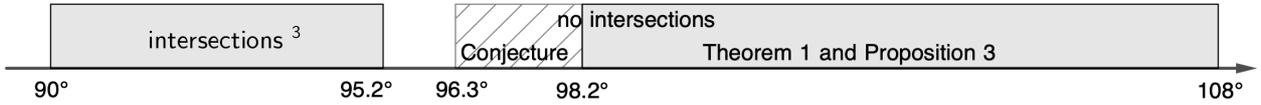

Figure 2: The goal of this article

## 2. Basics
### 2. 1. Defining the paperfolding process as a sequence of symbols

The paperfolding polygons are generated by iteratively folding a strip of paper in the same direction and analyzing the resulting creases. This process generates a sequence of up and downs or hills and valleys. Referring to the geometric interpretation (see Definition 3) a down (valley) is denoted as L and an up (hill) as R. If the strip is orientated so that the first crease forms a down (= L), the sequence of symbols starts with LLRLLRRLLLRRLRR … .

Definition 1
Let $\sigma: \mathbb{N} \to \{L, R\}$ be the function that assigns the symbol at position $n$ in the sequence.

This function allows us to define the paperfolding sequence.

Definition 2 (Inflation Law)
The paperfolding sequence consists of symbols L and R determined by:
$$\sigma(n) = \begin{cases} L & \text{if } n \equiv 1 \bmod 4 \\ R & \text{if } n \equiv 3 \bmod 4 \\ \sigma\left(\frac{n}{2}\right) & \text{if } n \equiv 0 \bmod 2 \end{cases} \quad n \in \mathbb{N}$$

There is an alternative method to construct the sequence.

Lemma 1 (Reflection Law)

The paperfolding sequence can be defined by
$\sigma(2^k) = L, k \in \mathbb{N}_0$ and $\sigma(2^k - d) = \overline{\sigma(2^k + d)},\ d \in \{1, 2, \ldots, 2^k - 1\}$,
where $\bar{S}$ denotes the opposite symbol to the symbol $S \in \{L, R\}$.

Proof
See R. Albers, Papierfalten [3], chapter 3

### 2. 2. Defining the paperfolding process as a sequence of polygons
For a more geometric perspective, consider the folded paper strip as a polygon where creases form vertices and paper segments form edges.

Definition 3 (sequence of paperfolding polygons)
The polygon is dependent on an angle $\alpha$, called the unfolding angle.
$Q_n$ is a polygon of $2^n + 1$ vertices $P_{n,0}$ to $P_{n,2^n}$ and their connecting segments $\overline{P_{n,k}P_{n,k+1}}, k \in \{0,1,2,\ldots,2^n-1\}$, meeting at the angle $\alpha$. (Figure 3)
$Q_0$ is the straight line from $P_{0,0}$ to $P_{0,1}$.

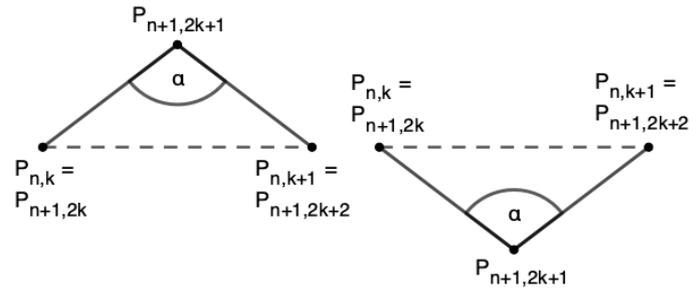

Figure 3: Geometric inflation

$Q_{n+1}$ is constructed from $Q_n$ by
  a) Renumbering the vertices $P_{n,k}$ to $P_{n+1,2k}$
  b) Replacing each segment $\overline{P_{n,k}P_{n,k+1}}$ by the sides of an isosceles triangle, meeting at $P_{n+1,2k+1}$ under the angle $\alpha$.
  c) The isosceles triangle is added
     i. to the left side, forming a right turn, if $\overline{P_{n,k}P_{n,k+1}}$ is for odd $k$ to even $k+1$.
     ii. to the right side, forming a left turn, if $\overline{P_{n,k}P_{n,k+1}}$ is for even $k$ to odd $k+1$.

Conclusion
  1. Each polygon begins at $P_{0,0}$ and ends at $P_{0,1}$.
  2. This geometric construction of $Q_n$ aligns with the inflation law (Definition 2).

Definition 4
  a) $\beta = |\sphericalangle P_{n,k+1}P_{n,k}P_{n+1,2k+1}|$
     $= \dfrac{180° - \alpha}{2} = 90° - \dfrac{\alpha}{2}$
  b) $q = \dfrac{|P_{n,k}P_{n+1,2k+1}|}{|P_{n,k}P_{n,k+1}|} = \dfrac{1}{2\cos\beta}$

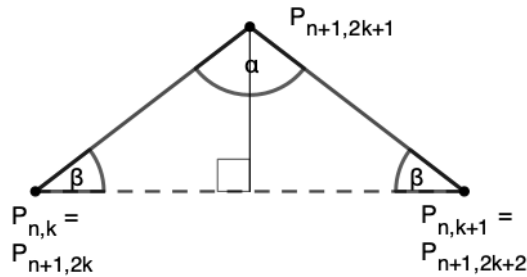

Figure 4: Definition of $\beta$ and $q$

With the parameters from Definition 4, two essential mappings can now be defined.

Definition 5
Mapping $\pi_0$: The rotation around $P_{0,0}$ by $-\beta$ combined with the dilation at center $P_{0,0}$ and scale factor $q$.
Mapping $\pi_1$: The rotation around $P_{0,0}$ by $\beta - 180°$ combined with the dilation at center $P_{0,0}$ and scale factor $q$, then the translation by $\overrightarrow{P_{0,0}P_{0,1}}$.

Remark
 Since the angles $-\beta$ and $\beta - 180°$ are both negative, the rotations are clockwise.

Lemma 2
The sequence of polygons $Q_n$ is generated recursively by
$Q_0 = \overline{P_{0,0}P_{0,1}}$ and $Q_{n+1} = \pi_0(Q_n) \cup \pi_1(Q_n)$.

Proof
For the special case $\beta = 45°$ see R. Albers, Papierfalten [3], chapter 8.

Remark
Lemma 2 corresponds to the reflection law (Lemma 1)

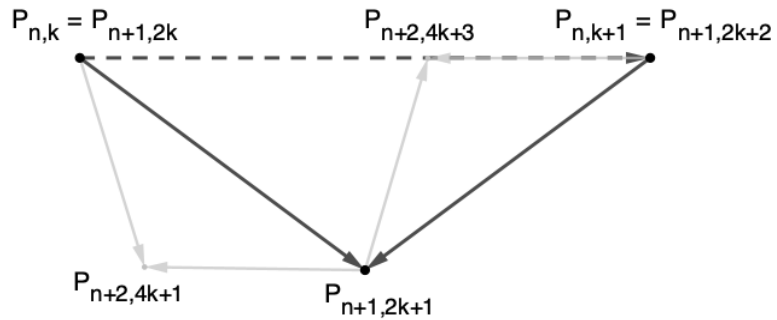

To facilitate further analysis, the segments get an orientation. The direction flows from vertices $P_{n,k}$ with an even index $k$ to those with an odd index $k+1$ or $k-1$. As illustrated in Figure 5, alternate segments are oriented retrograde to their predecessors.

Figure 5: Orientated segments ($k$ is even)

Lemma 3 (orientated segments)
When following the established segment orientation, all triangles in the construction of $Q_{n+1}$ from $Q_n$ are consistently inserted on the right-hand side of the directed segments.

Proof
This follows directly from Definition 3 c), see Figure 5

**Preview: Avoiding intersections of the paperfolding polygons $Q_n$**

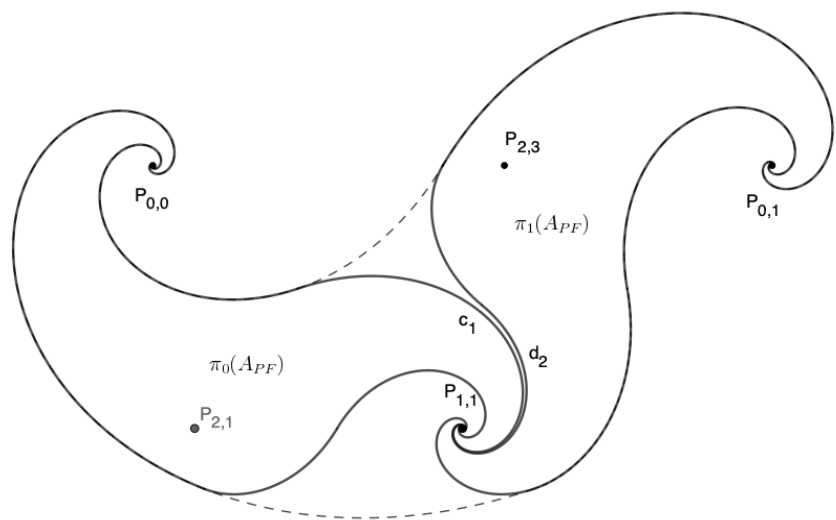

This article shows (in Proposition 3 and Theorem 1), that for unfolding angles $98.195° \leq \alpha \leq 108°$ the paperfolding polygons $Q_n$ avoid intersections.

The proof strategy relies on constructing a set $A_{PF}$ that contains all polygons $Q_n$ as subsets. If the images under the mappings $\pi_0$ and $\pi_1$ of this set are clearly separated,
i.e. $\pi_0(A_{PF}) \cap \pi_1(A_{PF}) = P_{1,1}$,
then no intersections occur in the paperfolding polygons.

Figure 6: The idea of the proof

The complete formal proof follows after Lemma 10.

# 3. Points on a straight line

Lemma 4a

Let $k$ be an odd integer. If there is a left turn at $P_{n,k}$, then the three points $P_{n,k-1}$, $P_{n,k} = P_{n+4,16k}$ and $P_{n+4,16k+1}$ are collinear. (Figure 7)

Proof

As $k$ is odd, the segment $\overline{P_{n,k}P_{n,k+1}}$ is orientated toward $P_{n,k}$. The inflation process, generating $Q_{n+1}$ from $Q_n$, adds $\overline{P_{n+1,2k}P_{n+1,2k+1}}$ so that it is rotated from $\overline{P_{n,k}P_{n,k+1}}$ by $\beta$. Subsequent inflation steps generate segments rotated by $-\beta$. So, the angle (absolute value) between $\overline{P_{n,k}P_{n,k+1}}$ and $\overline{P_{n+4,16k}P_{n+4,16k+1}}$ is $2\beta$.
As $\alpha + 2\beta = 180°$, the points $P_{n,k-1}$, $P_{n,k} = P_{n+4,16k}$ and $P_{n+4,16k+1}$ are on a straight line. $\square$

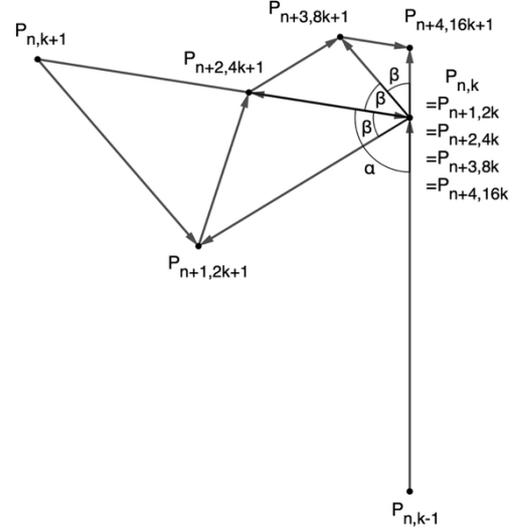

Figure 7: New segments after a left turn

Conclusion

The point $P_{n+4,16k+1}$ initiates another left turn (see Definition 2). By iterating this property, it yields $P_{n+8,256k+17}$ on the same straight line. This can be continued indefinitely.

Lemma 4b

Let $k$ be an odd integer. If there is a right turn at $P_{n,k}$, then the three points $P_{n,k+1}$, $P_{n,k} = P_{n+4,16k}$ and $P_{n+4,16k-1}$ are collinear. (Figure 8)

Proof

As $k$ is odd, the segment $\overline{P_{n,k}P_{n,k-1}}$ is orientated toward $P_{n,k}$. The inflation process, generating $Q_{n+1}$ from $Q_n$, adds $\overline{P_{n+1,2k}P_{n+1,2k-1}}$ so that it is rotated from $\overline{P_{n,k}P_{n,k-1}}$ $\beta$. Subsequent inflation steps generate segments rotated by $-\beta$. So, the angle (absolute value) between $\overline{P_{n,k}P_{n,k-1}}$ and $\overline{P_{n+4,16k}P_{n+4,16k-1}}$ is $2\beta$.
As $\alpha + 2\beta = 180°$, the points $P_{n,k+1}$, $P_{n,k} = P_{n+4,16k}$ and $P_{n+4,16k-1}$ are on a straight line. $\square$

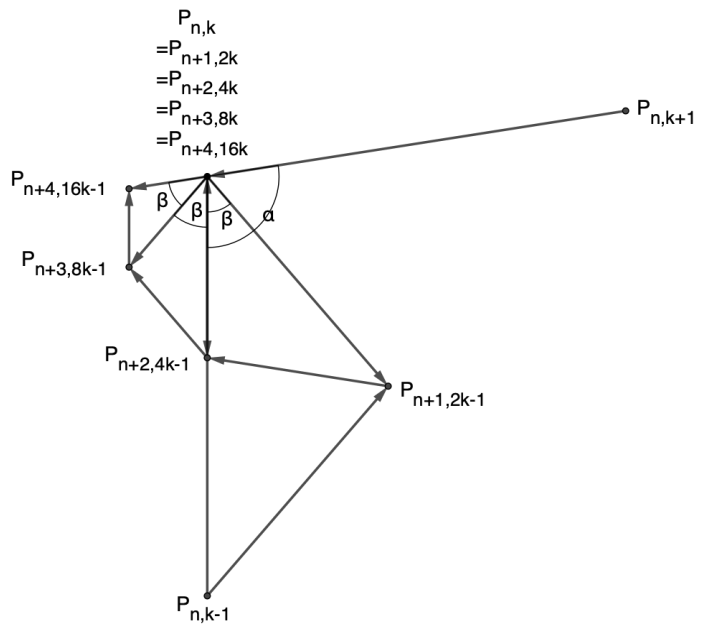

Figure 8: New segments before a right turn

Conclusion

The point $P_{n+4,16k-1}$ initiates another right turn (see Definition 2). By iterating this property, it yields $P_{n+8,256k-17}$ on the same straight line. This can be continued indefinitely.

In the following part of the article logarithmic spirals are intensively used. So, it is useful to denote some

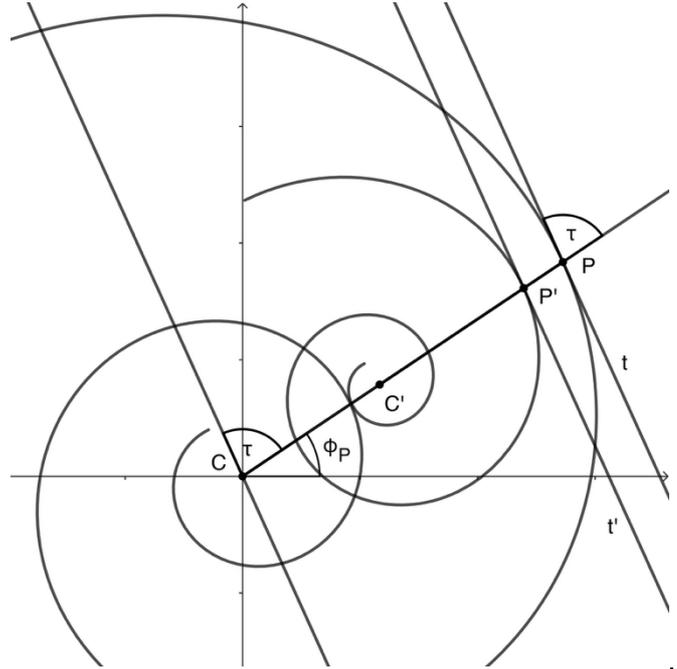

Basic facts for logarithmic spirals
In a polar coordinate system, a logarithmic spiral with the origin C as center is described by
$$r(\varphi) = ab^\varphi, a, b > 0$$
For every point P of the spiral the angle $\tau$ between the radius CP and the tangent to the spiral at P is constant and given by
$$\tan(\tau) = \frac{1}{\ln b}$$
This yields for two spirals with center C resp. C' and the same base $b$: If P and P' are two points of the spirals and collinear with C and C', then the tangents t and t' to the spirals at P and P' are parallel. If P = P', then the tangent is a common tangent to both spirals.

Let base $b > 1$ and $\varphi_P$ be the polar angle of P measured from center C. Let t be the tangent to the spiral at P. Then all points of the spiral with a polar angle $\varphi \leq \varphi_P + \tau$ are with respect to t in the half-plane containing center C.

The dilation of a spiral at its center C by the factor $s$ and the rotation around C by the angle $\rho$ transforms the describing formula from
$r(\varphi) = a\, b^\varphi$, $\varphi_{min} \leq \varphi \leq \varphi_{max}$ to $\tilde{r}(\varphi) = sa\, b^{\varphi - \rho}$, $\varphi_{min} + \rho \leq \varphi \leq \varphi_{max} + \rho$

## 4. Constructing a borderline for the polygons $Q_n$

Construction 1: The maximal distance to $P_{0,0}$
According to Lemma 4a) the three points $P_{n,k-1}$, $P_{n,k} = P_{n+4,16k}$ and $P_{n+4,16k+1}$ are collinear. After four inflation steps, this pattern repeats with additional points on the same line. If the length of $\overline{P_{n,k-1}P_{n,k}}$ is $L$, the total length of consecutive segments on this line follows a geometric series: $L(1 + q^4 + q^8 + q^{12} + \cdots) = L\frac{1}{1-q^4}$ (for $q$ see Definition 4).
Applying this to $\overline{P_{n,k-1}P_{n,k}} = \overline{P_{1,0}P_{1,1}}$ there exists a straight line from $P_{0,0}$ to $P_{1,1}$ to $P_{5,17}$, $P_{9,273}$, $P_{13,4369}$, ... (see Figure 9). Given $|P_{0,0}P_{0,1}| = 1 \Rightarrow |P_{1,0}P_{1,1}| = q$, the distance from $P_{0,0}$ through $P_{1,1}$ has in the limit the length of $q\frac{1}{1-q^4}$.
Then applying it to $\overline{P_{n,k-1}P_{n,k}} = \overline{P_{2,0}P_{2,1}}$ and compared to $\overline{P_{1,0}P_{1,1}}$, this segment is rotated by $-\beta$ (clockwise) and dilated by $q$. Therefore, the limit length becomes $q^2\frac{1}{1-q^4}$. Iterating this process generates a spiral, that serves as an outer boundary for the paperfolding polygons around $P_{0,0}$. In a polar system with origin in $P_{0,0}$ and $\overrightarrow{P_{0,0}P_{0,1}}$ defines the direction of $\varphi = 0$ this spiral is given by

$$a_0(\varphi) = \frac{1}{1-q^4} q^{-\frac{\varphi}{\beta}}, \varphi \leq -\beta$$

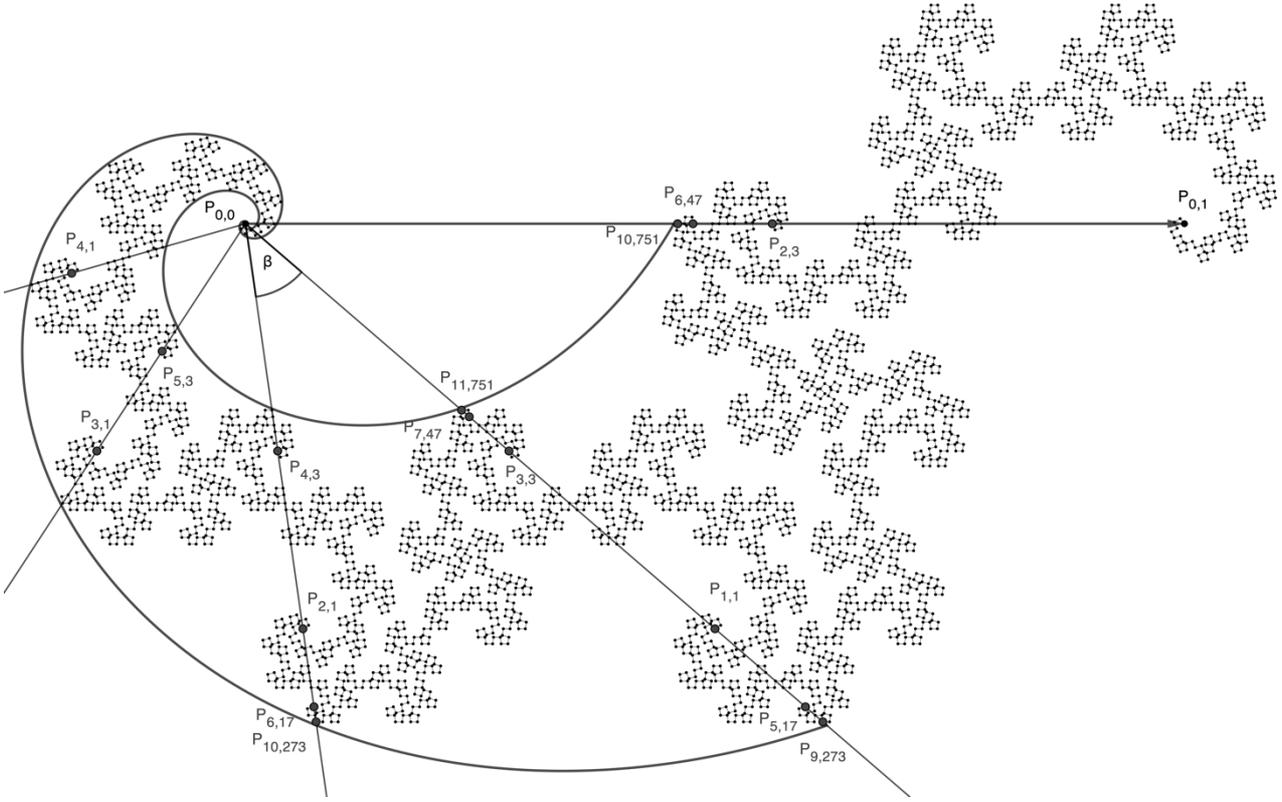

Figure 9: The polygon $Q_{12}$ with the outer and inner spiral at $P_{0,0}$

<u>Construction 2</u>: The minimal distance to $P_{0,0}$

According to Lemma 4b) the three points $P_{n,k+1}$, $P_{n,k} = P_{n+4,16k}$ and $P_{n+4,16k-1}$ are collinear. After four inflation steps, this pattern repeats with additional points on the same line. If the length of $\overline{P_{n,k+1}P_{n,k}}$ is $L$, the total length of consecutive segments on this line follows a geometric series: $L(1 + q^4 + q^8 + q^{12} + \cdots) = L\frac{1}{1-q^4}$.

Applying this to $\overline{P_{n,k+1}P_{n,k}} = \overline{P_{2,4}P_{2,3}}$ (remember $P_{2,4} = P_{0,1}$) there exists a straight line from $P_{0,1}$ to $P_{2,3}$ to $P_{6,47}$, $P_{10,751}$, ... (see Figure 9). Given $|P_{0,0}P_{0,1}| = 1 \Rightarrow |P_{2,3}P_{0,1}| = q^2$ and the length is in the limit $\frac{q^2}{1-q^4}$. As $P_{0,0}$ lies on the same line, the distance from this limit to $P_{0,0}$ is $\left(1 - \frac{q^2}{1-q^4}\right)$.

Applying this to the next right turn at $P_{n,k} = P_{3,3}$. It yields $P_{n,k+1} = P_{3,4} = P_{1,1}$. Relative to $\overline{P_{2,4}P_{2,3}}$ the straight line is rotated by $-\beta$ (clockwise) and everything is shortened by $q$, resulting in a minimal distance to $P_{0,0}$ of $q\left(1 - \frac{q^2}{1-q^4}\right)$. Iterating this process generates a spiral that serves as the minimum distance of the paperfolding polygons to $P_{0,0}$.

In the same coordinate system of construction 1 it is given by

$$b_0(\varphi) = \left(1 - \frac{q^2}{1-q^4}\right) q^{-\frac{\varphi}{\beta}}, \varphi \leq 0$$

Construction 3: The maximal and minimal distance to $P_{0,1}$

According to Lemma 2 every polygon $Q_n$ is generated by the mappings $\pi_0$ and $\pi_1$.
Since $P_{0,1} = \pi_1(P_{0,0})$, the spirals $a_0$ and $b_0$ are transformed by $\pi_1$.
$\pi_1$: dilation by $q$, rotation by $\beta - 180° = -\alpha - \beta$

$$a_0(\varphi) = \frac{1}{1-q^4} q^{-\frac{\varphi}{\beta}}, \varphi \le -\beta \xrightarrow{\pi_1} c_0(\varphi) = \frac{1}{1-q^4} q^{-\frac{\varphi+\alpha}{\beta}}, \varphi \le -180°$$

For the further argumentation the domain for $c_0$ is extended to $\varphi \le -180° + \beta = -\alpha - \beta$.

$$b_0(\varphi) = \left(1 - \frac{q^2}{1-q^4}\right) q^{-\frac{\varphi}{\beta}}, \varphi \le 0 \xrightarrow{\pi_1} d_0(\varphi) = \left(1 - \frac{q^2}{1-q^4}\right) q^{-\frac{\varphi+\alpha}{\beta}}, \varphi \le -\alpha - \beta$$

The origin of the underlying coordinate system for $c_0$ and $d_0$ is translated to $P_{0,1}$.

Construction 4: Closing the gap between $a_0$ and $d_0$

Since $P_{1,1} = \pi_1(P_{0,1})$, the spiral $c_0$ is transformed by $\pi_1$.
$\pi_1$: dilation by $q$, rotation by $\beta - 180° = -\alpha - \beta$

$$c_0(\varphi) = \frac{1}{1-q^4} q^{-\frac{\varphi+\alpha}{\beta}}, \varphi \le -180° \xrightarrow{\pi_1} e_0(\varphi) = \frac{1}{1-q^4} q^{-\frac{\varphi+2\alpha}{\beta}}, \varphi \le -360° + \beta$$

Applying the parameter-transformation $\varphi \to \varphi - 360° = \varphi - (2\alpha + 4\beta)$ yields

$$e_0(\varphi) = \frac{1}{1-q^4} q^{-\frac{\varphi-4\beta}{\beta}}, \varphi \le \beta$$

The origin of the underlying coordinate system is translated to $P_{1,1}$.

For a complete encompassing of the paperfolding polygons, only the part of the spiral is required until it reaches $a_0$ at the polar angle $-\beta$. Therefore, the interval for $\varphi$ is restricted to $-\beta \le \varphi \le \beta$.

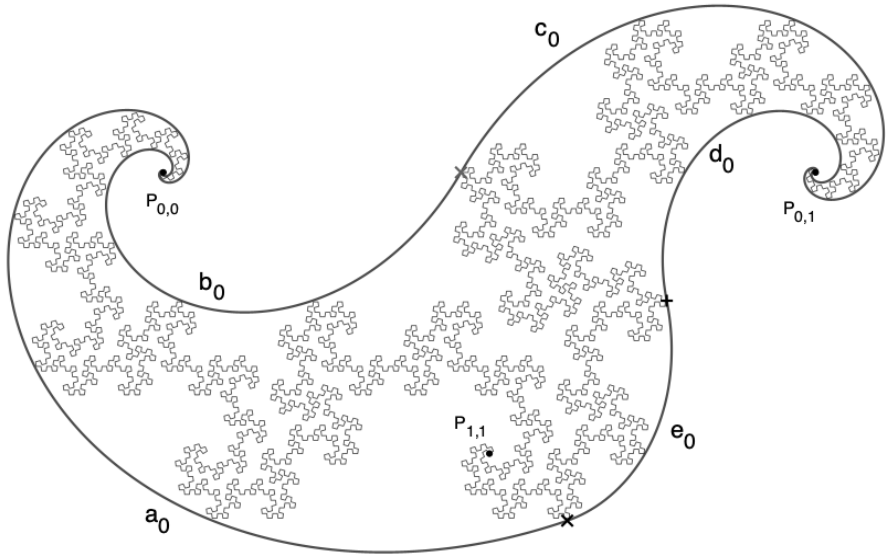

Figure 10: The five parts of the borderline

## 5. Defining a set, that is a hull for the paperfolding polygons

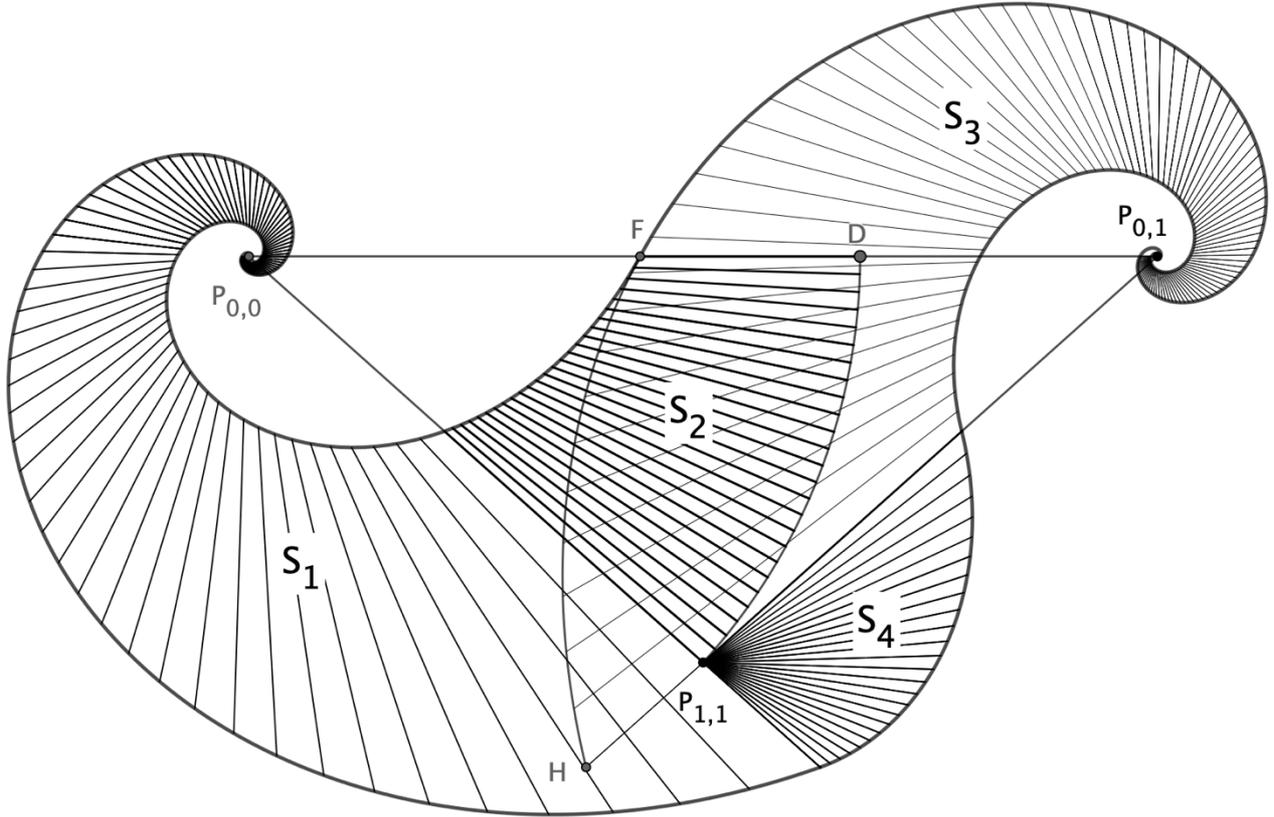

Figure 11: The four sets of points

Definition 6

a) F is the point on $\overline{P_{0,0}P_{0,1}}$ with $|P_{0,0}F| = b_0(0°) = 1 - \frac{q^2}{1-q^4}$.

D is the point on $\overline{P_{0,0}P_{0,1}}$ with $|P_{0,0}D| = q$.

H is the point on $c_0$ with polar angle $-\alpha - \beta$.

b) For $90° \leq \alpha \leq 108°$ (corresponding to $45° \geq \beta \geq 36°, \frac{1}{2}\sqrt{2} \geq q \geq \frac{\sqrt{5}-1}{2}$)
the set of points, that is delimited by the encompassing spirals constructed above, consist of four parts:

$(d_{i,k}; \varphi)$ are the polar coordinates of a point, $d_{i,k}$ is the distance to the point $P_{i,k}$

$S_1 = \{(d_{0,0}; \varphi) \in \mathbb{R}^2 | \ b_0(\varphi) \leq d_{0,0} \leq a_0(\varphi), \varphi \leq -\beta\}$
$S_2 = \{(d_{0,0}; \varphi) \in \mathbb{R}^2 | \ b_0(\varphi) \leq d_{0,0} \leq q, -\beta \leq \varphi \leq 0°\}$
$S_3 = \{(d_{0,1}; \varphi) \in \mathbb{R}^2 | \ d_0(\varphi) \leq d_{0,1} \leq c_0(\varphi), \varphi \leq -\alpha - \beta\}$
$S_4 = \{(d_{1,1}; \varphi) \in \mathbb{R}^2 | \ 0 \leq d_{1,1} \leq e_0(\varphi), -\beta \leq \varphi \leq \beta\}$

Then the total set of points is
$A_{PF} = S_1 \cup S_2 \cup S_3 \cup S_4$

This set constitutes a hull for all paperfolding polygons $Q_n$ of a fixed unfolding angle $\alpha$. Due to its characteristic shape, we will refer to it as the "duck-shaped hull" or, concisely, the "duck".

## Lemma 5

Let D and H be the points defined in Definition 6 (see Figure 11). For $\frac{1}{2}\sqrt{2} \geq q \geq \frac{\sqrt{5}-1}{2}$

a) $D \epsilon S_3$.

b) $H \epsilon S_1$

Proof

a) $D \epsilon S_3$ holds, if $d_0(-180°) \leq 1 - q \leq c_0(-180°)$

$\Leftrightarrow \left(1 - \frac{q^2}{1-q^4}\right) q^2 \leq 1 - q \leq \frac{1}{1-q^4} q^2$ which holds for $\frac{1}{2}\sqrt{2} \geq q \geq \frac{\sqrt{5}-1}{2}$.

b) $H \epsilon S_1$ holds, if $a_0(\lambda) \leq |P_{0,0}H| \leq b_0(\lambda)$, where $\lambda$ is the polar angle for $\overline{P_{0,0}H}$.

H is defined by $|P_{0,1}H| = c_0(-\alpha - \beta)$. Applying the law of cosine in the triangle $P_{0,0}HP_{0,1}$

yields $|P_{0,0}H| = q \frac{\sqrt{1 - q^2 + q^6}}{1-q^4}$ and $\tan|\lambda| = \frac{\sqrt{4q^2-1}}{1-2q^4}$. A numerical analysis shows that

$a_0\left(-\arctan\left(\frac{\sqrt{4q^2-1}}{1-2q^4}\right)\right) \leq q \frac{\sqrt{1 - q^2 + q^6}}{1-q^4} \leq b_0\left(-\arctan\left(\frac{\sqrt{4q^2-1}}{1-2q^4}\right)\right)$ holds for $\frac{1}{2}\sqrt{2} \geq q \geq \frac{\sqrt{5}-1}{2}$.

## Definition 7

Generalization of $A_{PF}$

Let $\pi$ be the composition of the maps $\pi_0$ and $\pi_1$ such that

$\overline{P_{n,k}P_{n,k+1}} = \pi(\overline{P_{0,0}P_{0,1}})$, $k \in \{0,1,2, \ldots, 2^n - 1\}$.

Then $A_{PF}(\overline{P_{n,k}P_{n,k+1}})$ is defined by $A_{PF}(\overline{P_{n,k}P_{n,k+1}}) = \pi(A_{PF})$.

So, $A_{PF}(\overline{P_{n,k}P_{n,k+1}})$ is a contracted copy of $A_{PF}$ orientated at the segment $\overline{P_{n,k}P_{n,k+1}}$. The initial set $A_{PF}$ can be written as $A_{PF}(\overline{P_{0,0}P_{0,1}})$.

For further argumentation additional subsets of $A_{PF}$ are defined.

## Definition 8 (see Figure 12)

a) Points:

A is the intersection of the straight line $P_{0,0}P_{1,1}$ with $a_0$, $|P_{0,0}A| = a_0(-\beta)$.

E is the intersection of the straight line $P_{0,0}P_{1,1}$ with $a_0$, $|P_{0,0}E| = a_0(-180° - \beta)$.

B is the intersection of the straight line $P_{0,0}P_{1,1}$ with $b_0$, $|P_{0,0}B| = b_0(-\beta)$.

C is the intersection of the straight line $P_{0,1}P_{1,1}$ with $d_0$, $|P_{0,1}C| = d_0(-\alpha - \beta)$.

b) Lines:

As $|P_{1,1}B| = |P_{1,1}C| = \frac{q^3}{1-q^4}$ there exists a circular arc $\widehat{BC}$ with midpoint $P_{1,1}$.

t is the tangent to $b_0$ at B.

$e_0'$ extends $e_0$ to cover polar angles from $-\beta$ to $-180° - \beta$ (center $P_{1,1}$).

The intersection of arc $\widehat{BC}$ with the tangent t is G.

c) Sets:

$T_1 = \{(d_{1,1}; \varphi) \in \mathbb{R}^2 | 0 \leq d_{1,1} \leq e_0'(\varphi), -\beta \geq \varphi \geq -180° - \beta \}$

$T_2$ is the set of points bounded by arc $\widehat{CG}$ (center $P_{1,1}$, radius $\frac{q^3}{1-q^4}$), segment $\overline{GB}$, segment $\overline{BP_{1,1}}$ and segment $\overline{P_{1,1}C}$. Therefore, every point of $T_2$ has a polar angle $\varphi$ to $P_{1,1}$ which satisfies $-180° - \beta \geq \varphi \geq -180° - \beta - \alpha$.

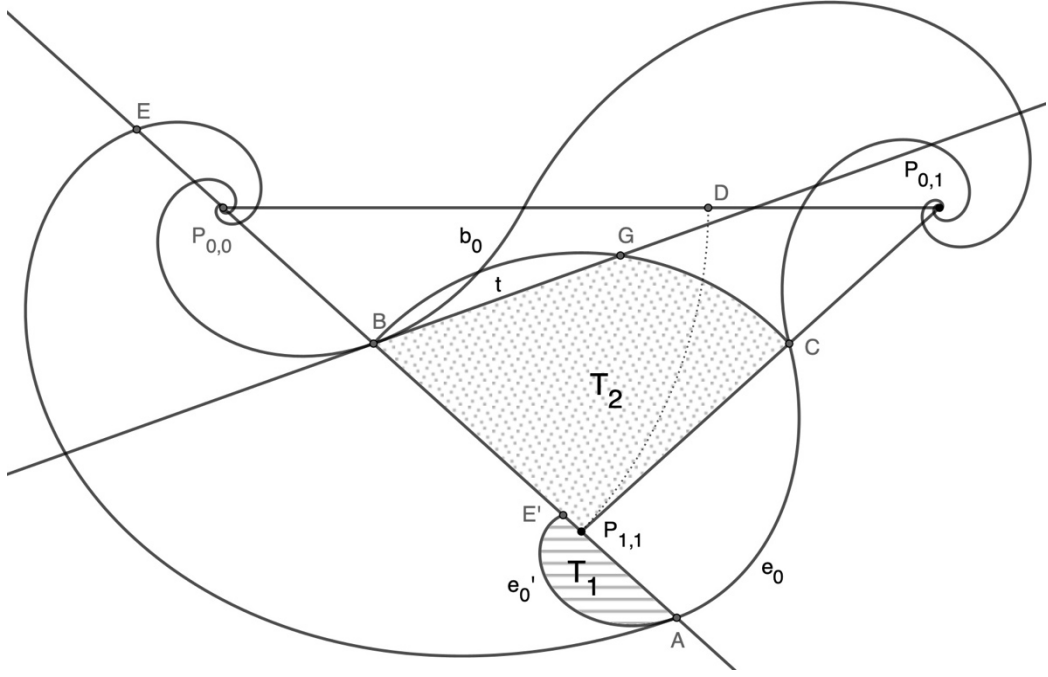

Figure 12: Subsets of $A_{PF}$

## Lemma 6

For $90° \leq \alpha \leq 108°$, $\frac{1}{2}\sqrt{2} \geq q \geq \frac{\sqrt{5}-1}{2}$ (see Figure 13)
  a) arc $\widehat{BC} \cap \overline{P_{0,0}P_{0,1}} = \emptyset$
  b) Point G is an inner point on arc $\widehat{BC}$.

## Proof

  a) The radius of arc $\widehat{BC}$ is $\frac{q^3}{1-q^4}$, the distance of $P_{1,1}$ to $\overline{P_{0,0}P_{0,1}}$ is $q \sin \beta = q\sqrt{1 - \frac{1}{4q^2}}$.
  
  The condition $\frac{q^3}{1-q^4} < q\sqrt{1 - \frac{1}{4q^2}}$ can be simplified to
  
  $1 - 4q^2 - 2q^4 + 12q^6 + q^8 - 4q^{10} < 0$, which holds for $0.724 \geq q \geq 0.524$.

  b) $\tau$ is the angle between the tangent to the spiral and the ray $\overrightarrow{P_{0,0}B}$.
  Therefore, $\tan \tau = -\frac{\beta}{\ln q} = \frac{\beta}{\ln(2\cos\beta)}$
  G lies between C and B is equivalent to $\tau > \beta$.
  This can be transformed to
  $\frac{\beta}{\ln(2\cos\beta)} > \tan \beta$.
  This holds for all $0° < \beta < 45°$

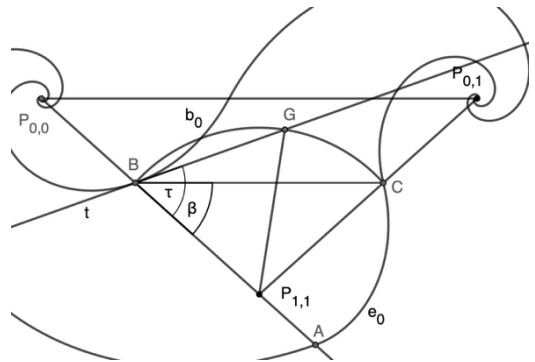

Figure 13: Angle $\tau$ and point G

## Lemma 7
  a) $T_1 \subset S_1 \subset A_{PF}$
  b) $T_2 \subset (S_2 \cup S_3) \subset A_{PF}$

Proof (see Figure 12)

a) Let $T_1^*$ be the area bounded by the spiral arc $\widehat{AE}$ and the line segment $\overline{EA}$. The dilation with center A and scaling factor $q^4$ maps $T_1^*$ to $T_1$ because:

$$|AP_{0,0}| = a_0(-\beta) = \frac{1}{1-q^4} q^{-\frac{-\beta}{\beta}} = \frac{1}{1-q^4} q \quad |AP_{1,1}| = e_0(-\beta) = \frac{1}{1-q^4} q^{-\frac{-\beta-4\beta}{\beta}} = \frac{1}{1-q^4} q^5$$

So, the center $P_{0,0}$ is mapped to the center $P_{1,1}$.

$$|P_{0,0}P| = a_0(\varphi) = \frac{1}{1-q^4} q^{-\frac{\varphi}{\beta}} \quad |P_{1,1}P'| = e_0'(\varphi) = \frac{1}{1-q^4} q^{-\frac{\varphi-4\beta}{\beta}} = \frac{1}{1-q^4} q^{-\frac{\varphi}{\beta}} q^4$$

So, every point P on $a_0$ is mapped to a point P' on $e_0'$.
$T_1^*$ is convex and $A \epsilon T_1^*$, so $T_1 \subset T_1^*$.
Now it has to be shown for all points of $T_1$ that the distance to $P_{0,0}$ is larger than the distance of an according point on $b_0$ to $P_{0,0}$.
Let P be a point of $T_1$. Applying the triangle inequality:

$$|P_{0,0}P| \geq |P_{0,0}P_{1,1}| - |P_{1,1}P| \geq q - |P_{1,1}A| = q - \frac{1}{1-q^4} q^5$$

$$\max(b_0(\varphi): -\beta \geq \varphi \geq -180° - \beta) = b_0(-\beta) = \left(1 - \frac{q^2}{1-q^4}\right) q = q - \frac{1}{1-q^4} q^3$$

Therefore, $|P_{0,0}P| \geq \max(b_0(\varphi): -\beta \geq \varphi \geq -180° - \beta)$, which implies $T_1 \subset S_1$. □

b) Let $P \epsilon T_2$
Case 1: $|P_{0,0}P| \leq q$
The definition of $S_2$ directly yields $P \epsilon S_2$.
Case 2: $|P_{0,0}P| > q$
The definition of $T_2$ yields $|P_{1,1}P| \leq \frac{q^3}{1-q^4}$.
Now it has to be shown for P that the distance to $P_{0,1}$ is equal or larger than the distance of an according point on $d_0$ to $P_{0,1}$. Applying the triangle inequality:

$$|PP_{0,1}| \geq |P_{1,1}P_{0,1}| - |P_{1,1}P| \geq q - \frac{q^3}{1-q^4}$$

$$\max (d_0(\varphi): -180° + \beta \geq \varphi \geq -180°) = d_0(-180° + \beta) = q - \frac{q^3}{1-q^4}$$

Therefore, $|PP_{0,1}| \geq \max (d_0(\varphi): -180° + \beta \geq \varphi \geq -180°)$ which implies $P \epsilon S_3$. □

Proposition 1
a) $\pi_0(A_{PF}) \subset A_{PF}$
b) $\pi_1(A_{PF}) \subset A_{PF}$

Proof a)
Since $\pi_0$ preserves the ratio, it suffices to examine the mapping of the encompassing spirals.
$\pi_0$: A dilation by $q$ and a rotation by $-\beta$, both with center $P_{0,0}$
points:   $P_{0,0} \to P_{0,0}$   $P_{0,1} \to P_{1,1}$   $P_{1,1} \to P_{2,1}$
spirals:

$a_0(\varphi) = \frac{1}{1-q^4} q^{-\frac{\varphi}{\beta}}, \varphi \leq -\beta$, center $P_{0,0}$

$\to a_1(\varphi) = \frac{1}{1-q^4} q^{-\frac{\varphi}{\beta}}, \varphi \leq -2\beta$, center $P_{0,0}$

$b_0(\varphi) = \left(1 - \frac{q^2}{1-q^4}\right) q^{-\frac{\varphi}{\beta}}, \varphi \leq 0$, center $P_{0,0}$

$\to b_1(\varphi) = \left(1 - \frac{q^2}{1-q^4}\right) q^{-\frac{\varphi}{\beta}}, \varphi \leq -\beta$, center $P_{0,0}$

$$c_0(\varphi) = \frac{1}{1-q^4} q^{-\frac{\varphi+\alpha}{\beta}}, \varphi \leq -180°, \text{center } P_{0,1}$$

$$\rightarrow c_1(\varphi) = \frac{1}{1-q^4} q^{-\frac{\varphi+\alpha}{\beta}}, \varphi \leq -180° - \beta, \text{center } P_{1,1}$$

$$d_0(\varphi) = \left(1 - \frac{q^2}{1-q^4}\right) q^{-\frac{\varphi+\alpha}{\beta}}, \varphi \leq -\alpha - \beta, \text{center } P_{0,1}$$

$$\rightarrow d_1(\varphi) = \left(1 - \frac{q^2}{1-q^4}\right) q^{-\frac{\varphi+\alpha}{\beta}}, \varphi \leq -180°, \text{center } P_{1,1}$$

$$e_0(\varphi) = \frac{1}{1-q^4} q^{-\frac{\varphi-4\beta}{\beta}}, \quad -\beta \leq \varphi \leq \beta, \text{center } P_{1,1}$$

$$\rightarrow e_1(\varphi) = \frac{1}{1-q^4} q^{-\frac{\varphi-4\beta}{\beta}}, \quad -2\beta \leq \varphi \leq 0°, \text{center } P_{2,1}$$

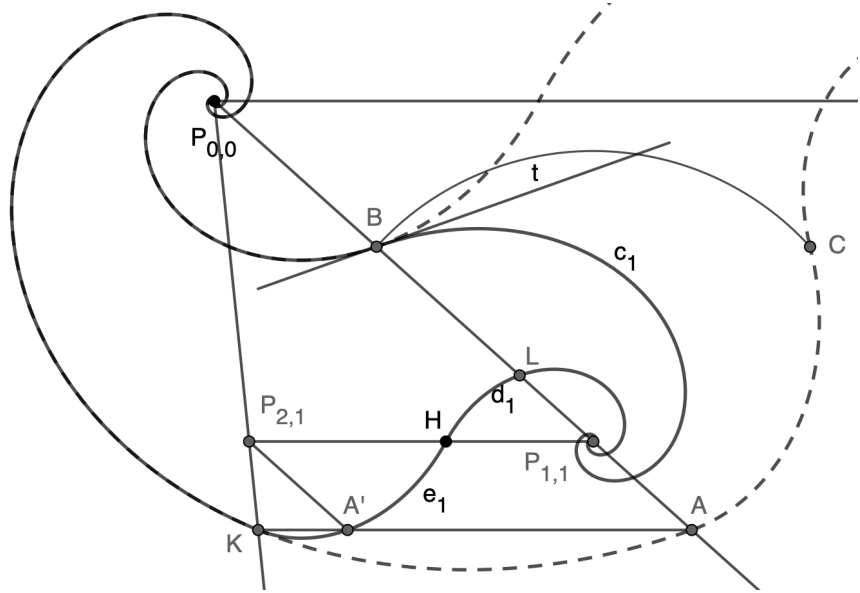

Figure 14: $\pi_0(A_{PF})$

The dilation by $q$ and the rotation by $-\beta$ cancel each other out. The spirals are just translated to different centers and the interval for $\varphi$ is modified.
$a_0$ and $b_0$ are mapped to themselves.

$c_1$ is tangent to $b_0$ at B with common tangent line t. $|P_{0,0}B| = b_0(-\beta) = \left(1 - \frac{q^2}{1-q^4}\right) q$

$|BP_{1,1}| = c_1(-180° - \beta) = \frac{1}{1-q^4} q^{-\frac{-180°-\beta+\alpha}{\beta}} = \frac{1}{1-q^4} q^3$, yielding $|P_{0,0}B| + |BP_{1,1}| = q$.

For $-180° - \beta \geq \varphi \geq -180° - \beta - \alpha$   $c_1$ is a subset of $T_2$, as
$c_1(\varphi) \leq |BP_{1,1}| = \frac{1}{1-q^4} q^3$ and $c_1$ is in the same half-plane as $P_{1,1}$ with respect to t.

For $-180° - \beta - \alpha = -360° + \beta \geq \varphi \geq -360° - \beta$   $c_1$ is a subset of $S_4$ as
$c_1(\varphi) = \frac{1}{1-q^4} q^{-\frac{\varphi+\alpha}{\beta}} \leq e_0(\varphi + 360°) = \frac{1}{1-q^4} q^{-\frac{\varphi+360°-4\beta}{\beta}} = \frac{1}{1-q^4} q^{-\frac{\varphi+2\alpha}{\beta}}$

For $-360° - \beta \geq \varphi \geq -540° - \beta$   $c_1$ is a subset of $T_1$ as
$c_1(\varphi) = \frac{1}{1-q^4} q^{-\frac{\varphi+\alpha}{\beta}} \leq e_0'(\varphi + 360°) = \frac{1}{1-q^4} q^{-\frac{\varphi+2\alpha}{\beta}}$

For $-540° - \beta \geq \varphi$ the argumentation repeats with contracted radii.

For $-180° \geq \varphi \geq -180° - \beta$  $d_1$ is a subset of $S_1$ as:
$d_1$ lies in the same half-plane as $P_{0,0}$ with respect to $P_{2,1}P_{1,1}$, so only the distance to $b_0$ needs to be tested. Let P be a point on $d_1$. Applying the triangle inequality:
$$|PP_{0,0}| \geq |P_{1,1}P_{0,0}| - |P_{1,1}P| \geq q - d_1(-180°) = q - \left(1 - \frac{q^2}{1-q^4}\right)q^{\frac{180°-\alpha}{\beta}} = q - \left(1 - \frac{q^2}{1-q^4}\right)q^2$$
The polar angle of line $PP_{0,0}$ is between $-\beta$ and $-2\beta$.
In this domain $b_0(\varphi) \leq \left(1 - \frac{q^2}{1-q^4}\right)q$. Numerical analysis implies that
$b_0(\varphi) \leq \left(1 - \frac{q^2}{1-q^4}\right)q \leq q - \left(1 - \frac{q^2}{1-q^4}\right)q^2 \leq |PP_{0,0}|$ holds for $q > 0.57$, corresponding to
$\alpha \leq 122°$.
For $-180° - \beta \geq \varphi$  $d_1$ is a subset of $T_2 \cup S_4 \cup T_1$ as
$$d_1(\varphi) = \left(1 - \frac{q^2}{1-q^4}\right)q^{-\frac{\varphi+\alpha}{\beta}} \leq c_1(\varphi) = \frac{1}{1-q^4}q^{-\frac{\varphi+\alpha}{\beta}}$$

$e_1$:
Let K and A be the intersection points of $P_{0,0}P_{2,1}$ and $P_{0,0}P_{1,1}$ with $a_0$. The straight-line KA intersects $e_1$ in A'. The dilation with center K and factor $q^4$ maps the spiral arc $\widehat{KA}$ (part of $a_0$) to $\widehat{KA'}$ (part of $e_1$):
$$|KP_{0,0}| = a_0(-2\beta) = \frac{1}{1-q^4}q^2 \xrightarrow{\text{dilation}} |KP_{0,0}|q^4 = \frac{1}{1-q^4}q^6 = e_1(-2\beta) = |KP_{2,1}|$$
So, $P_{0,0}$, the center of $a_0$ is mapped to $P_{2,1}$, the center of $e_1$.
$$a_0(\varphi) = \frac{1}{1-q^4}q^{-\frac{\varphi}{\beta}} \xrightarrow{\text{dilation}} \frac{1}{1-q^4}q^{-\frac{\varphi}{\beta}}q^4 = \frac{1}{1-q^4}q^{-\frac{\varphi-4\beta}{\beta}} = e_1(\varphi)$$
So, every point on the spiral arc $\widehat{KA}$ is mapped to a point on $\widehat{KA'}$.
Let $T_4$ be the area, bounded by the arc $\widehat{KA}$ and the line segment $\overline{KA}$. $T_4$ is convex and is a subset of $S_1$. Since $K \epsilon T_4$ the image arc $\widehat{KA'}$ is a subset of $T_4$.
So, $\widehat{KA'} = \left\{(d_{2,1}; \varphi) \in \mathbb{R}^2 \mid d_{2,1} = \frac{1}{1-q^4}q^{-\frac{\varphi+4\beta}{\beta}}, -2\beta \geq \varphi \geq -\beta\right\}$
Let P be a point on the remaining spiral arc $\widehat{HA'}$. Then
  i. The polar angle $\varphi$ of P: $-\beta \leq \varphi \leq 0°$.
  ii. The distance to $P_{2,1}$: $|P_{2,1}P| \leq |P_{2,1}H| = e_1(0°) = \frac{1}{1-q^4}q^4 < q^2 = |P_{2,1}P_{1,1}|$,
      which holds for $q \leq \frac{1}{2}\sqrt{2}$.
  iii. P lies with respect to the tangent to $e_1$ at A in the halfplane that contains $P_{2,1}$.
The three properties yield that the spiral arc $\widehat{HA'}$ is subset of the parallelogram $P_{1,1}P_{2,1}A'A \subset S_1$ (for $\overline{P_{2,1}P_{1,1}} \subset S_1$  see Lemma 9).

Therefore, $\pi_0(A_{PF}) \subset A_{PF}$. □

Proof b)
Since $\pi_1$ preserves the ratio, it suffices to examine the mapping of the encompassing spirals.
$\pi_1$: A dilation by q and a rotation by $\beta - 180°$, both with center $P_{0,0}$, and a translation by $\overrightarrow{P_{0,0}P_{0,1}}$.
points:   $P_{0,0} \to P_{0,1}$   $P_{0,1} \to P_{1,1}$   $P_{1,1} \to P_{2,3}$

spirals:

$a_0(\varphi) = \dfrac{1}{1-q^4} q^{-\frac{\varphi}{\beta}}, \varphi \leq -\beta$ , center $P_{0,0}$

$\rightarrow \quad a_2(\varphi) = \dfrac{1}{1-q^4} q^{-\frac{\varphi+\alpha}{\beta}}, \varphi \leq -180°$ , center $P_{0,1}$

$b_0(\varphi) = \left(1 - \dfrac{q^2}{1-q^4}\right) q^{-\frac{\varphi}{\beta}}, \varphi \leq 0$ , center $P_{0,0}$

$\rightarrow \quad b_2(\varphi) = \left(1 - \dfrac{q^2}{1-q^4}\right) q^{-\frac{\varphi+\alpha}{\beta}}, \varphi \leq -\alpha-\beta$ , center $P_{0,1}$

$c_0(\varphi) = \dfrac{1}{1-q^4} q^{-\frac{\varphi+\alpha}{\beta}}, \varphi \leq -180°$ , center $P_{0,1}$

$\rightarrow \quad c_2(\varphi) = \dfrac{1}{1-q^4} q^{-\frac{\varphi+2\alpha}{\beta}}, \varphi \leq \beta - 360°$ , center $P_{1,1}$

$d_0(\varphi) = \left(1 - \dfrac{q^2}{1-q^4}\right) q^{-\frac{\varphi+\alpha}{\beta}}, \varphi \leq -\alpha-\beta$ , center $P_{0,1}$

$\rightarrow \quad d_2(\varphi) = \left(1 - \dfrac{q^2}{1-q^4}\right) q^{-\frac{\varphi+2\alpha}{\beta}}, \varphi \leq -\alpha-180°$ , center $P_{1,1}$

$e_0(\varphi) = \dfrac{1}{1-q^4} q^{-\frac{\varphi-4\beta}{\beta}}, \quad -\beta \leq \varphi \leq \beta$ , center $P_{1,1}$

$\rightarrow \quad e_2(\varphi) = \dfrac{1}{1-q^4} q^{-\frac{\varphi+\alpha-4\beta}{\beta}}, \quad -180° \leq \varphi \leq -\alpha$ , center $P_{2,3}$

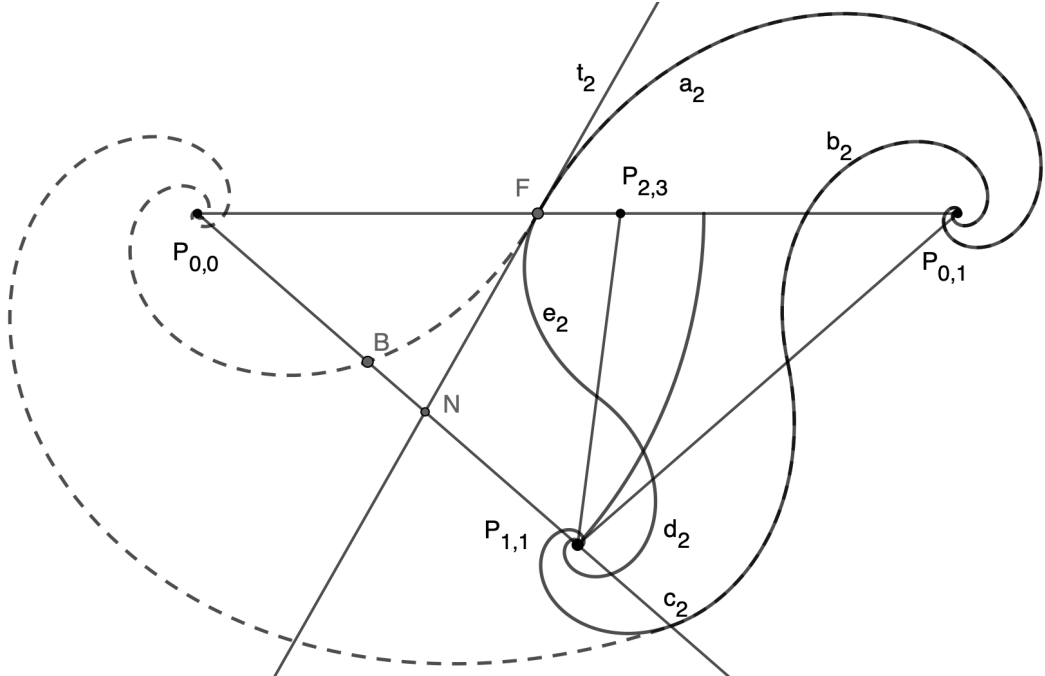

Figure 15: $\pi_1(A_{PF})$

$a_2$ is mapped to $c_0$ , $b_2$ is mapped to $d_0$ (see construction 3).

For $\beta \geq \varphi \geq -\beta$ , $c_2$ is mapped to $e_0$ (see construction 4) with $c_2(\varphi) = e_0(\varphi + 360°)$.
For $-\beta \geq \varphi \geq -180° - \beta \quad c_2$ is a subset of $T_1$
For $-180° - \beta \geq \varphi \quad c_2$ is a subset of $T_2 \cup S_4 \cup T_1$ (analogous to $c_1$ in Part a)).

For $-180° - \alpha \geq \varphi \geq -360° + \beta$, $d_2$ is a subset of $T_2$ as
$d_2(\varphi) \leq d_2(-180° - \alpha) = \left(1 - \frac{q^2}{1-q^4}\right) q^2 \leq \frac{q^3}{1-q^4}$.
Numerical analysis implies that this holds for $q > 0.57$ corresponding to $\alpha \leq 122°$.
For $-360° + \beta \geq \varphi$  $d_2$ is a subset of $T_2 \cup S_4 \cup T_1$, as $d_2(\varphi) \leq c_2(\varphi)$.

$e_2$ is tangent to $b_0$ at F, as $|P_{0,0}F| = b_0(0°) = 1 - \frac{q^2}{1-q^4}$
$|FP_{2,3}| = e_2(-180°) = q \frac{1}{1-q^4} q^5$, $|P_{2,3}P_{0,1}| = q^2$ and $|P_{0,0}F| + |FP_{2,3}| + |P_{2,3}P_{0,1}| = 1$.
All points of $e_2$ lie within the quadrilateral $P_{1,1}P_{2,3}FN$ which is for $90° \leq \alpha \leq 108°$ a subset of $S_2$.

Therefore, $\pi_1(A_{PF}) \subset A_{PF}$. □

Both maps together yield
$$\pi_0(A_{PF}) \cup \pi_1(A_{PF}) \subset A_{PF}$$

Corollary 1 (for generalized $A_{PF}$)

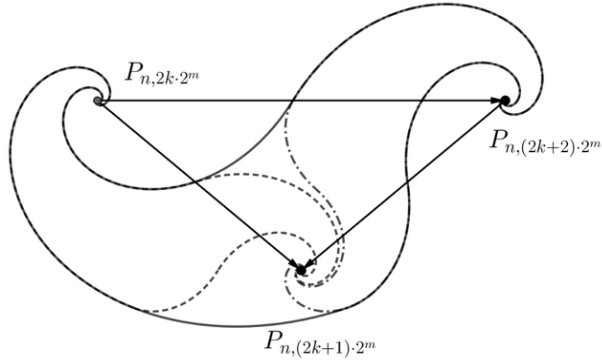

Figure 16: The mapping of generalized $A_{PF}$ by $\pi_0$ and $\pi_1$

$A_{PF}\left(\overline{P_{n,2k\cdot 2^m}P_{n,(2k+1)\cdot 2^m}}\right) \subset A_{PF}\left(\overline{P_{n,2k\cdot 2^m}P_{n,(2k+2)\cdot 2^m}}\right)$
$A_{PF}\left(\overline{P_{n,(2k+2)\cdot 2^m}P_{n,(2k+1)\cdot 2^m}}\right) \subset A_{PF}\left(\overline{P_{n,2k\cdot 2^m}P_{n,(2k+2)\cdot 2^m}}\right)$
$k, m, n \in \mathbb{N}_0, n \geq 1, 0 \leq m \leq n - 1, 0 \leq k \leq 2^{n-m-1} - 1$

## 6. The paperfolding polygons $Q_n$ are subsets of $A_{PF}$

Or simply said: All dragons are inside the duck.

Proposition 2
$$Q_n \backslash \left(\overline{P_{n,0}P_{n,1}} \cup \overline{P_{n,2^n-1}P_{n,2^n}}\right) \subset A_{PF}, n \in \mathbb{N}_0$$

The proof consists of two parts. First, all vertices of a paperfolding polygon are elements of $A_{PF}$ (Lemma 8) and second, all segments excluding the first $\overline{P_{n,0}P_{n,1}}$ and the last $\overline{P_{n,2^n-1}P_{n,2^n}}$ are subsets of $A_{PF}$ (Lemma 10).

Lemma 8
For all paperfolding polygons $Q_n, n \in \mathbb{N}_0$ the vertices $P_{n,k}, k \in \{0,1,2,\ldots,2^n\}$ are elements of $A_{PF}$.

Proof (by mathematical induction)
Base step: $P_{0,0}, P_{0,1}$ and $P_{1,1}$ are elements of $A_{PF}$.
Induction step: Let $P_{n,k}$ be an element of $A_{PF}$.
Then $\pi_0(P_{n,k}) \in \pi_0(A_{PF}) \subset A_{PF}$ and $\pi_1(P_{n,k}) \in \pi_1(A_{PF}) \subset A_{PF}$.
All vertices of the paperfolding polygons $Q_n, n \geq 2$ are obtained through iterative applications of $\pi_0$ and $\pi_1$ to $P_{1,1}$. □

Regarding the segments: Those connected to $P_{1,1}$ are special, since they are pictures under $\pi_0$ or $\pi_1$ of segments, that end in $P_{0,1}$ and theses segments are not subsets of $A_{PF}$.

Lemma 9
For $90° \leq \alpha \leq 108°$ all segments $\overline{P_{n,2^{n-1}-1}P_{n,2^{n-1}}}$ and $\overline{P_{n,2^{n-1}}P_{n,2^{n-1}+1}}, n \geq 2$ are subsets of $A_{PF}$. These are all segments which start or end at $P_{1,1} = P_{n,2^{n-1}}$.

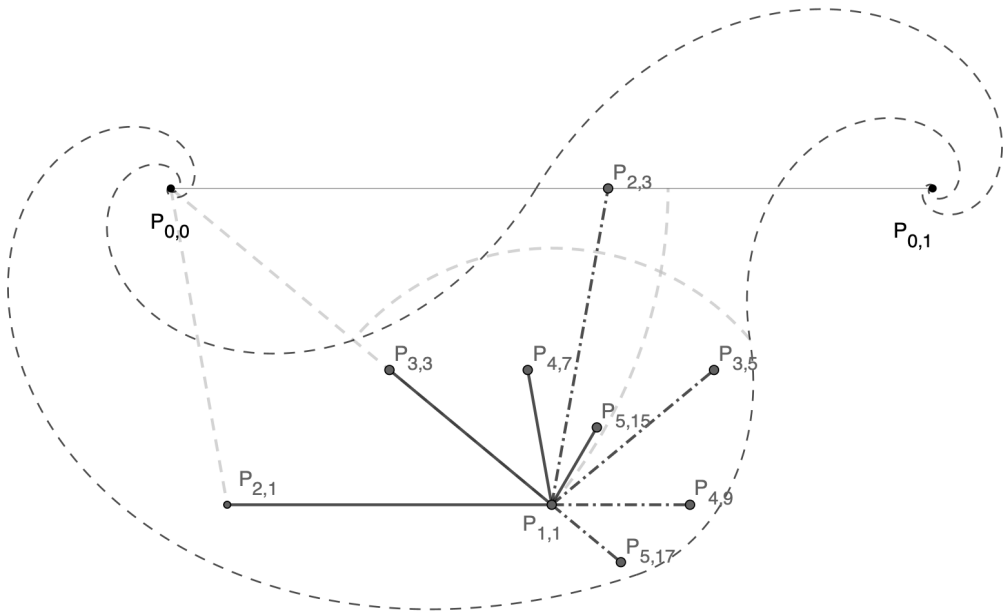

Figure 17: Some segments connected to P1,1

Proof for $\overline{P_{n,2^{n-1}-1}P_{n,2^{n-1}}}, n \geq 2$
With respect to $P_{1,1}$ the points $P_{n,2^{n-1}-1}$ have polar angle of $\varphi_n = -180° - (n-2)\beta$ and distance $r_n = q^n$.

**n = 2** Proof for $\overline{P_{2,1}P_{2,2}}$

Let P be a point on $\overline{P_{2,1}P_{2,2}}$. By the law of sines in triangle $PP_{2,2}P_{0,0}$ (Figure 18)

$$|P_{0,0}P| = s(\varphi) = \frac{q \sin(\beta)}{\sin(180° - |\varphi|)} = \frac{q \sin(\beta)}{-\sin(\varphi)}$$

Segment $\overline{P_{2,1}P_{2,2}}$ is a subset of $A_{PF}$,
if $b_0(\varphi) < s(\varphi) < a_0(\varphi)$ for $-2\beta \leq \varphi \leq -\beta$
$s(-\beta) = q < a_0(-\beta) = \frac{1}{1-q^4} q$
$s(-2\beta) = q^2 < a_0(-2\beta) = \frac{1}{1-q^4} q^2$

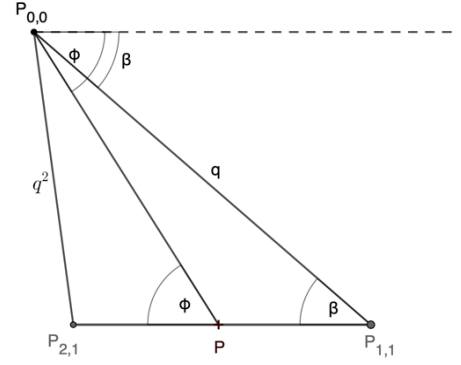

Figure 18: The segment $\overline{P_{2,1}P_{2,2}}$ (note: $P_{2,2} = P_{1,1}$)

Since $a_0(\varphi)$ forms a convex subset of $A_{PF}$ (with respect to $P_{0,0}$), $s(\varphi) < a_0(\varphi)$ for $-2\beta \leq \varphi \leq -\beta$.

For $-2\beta \leq \varphi \leq -\beta$ there is $max\{b_0(\varphi)\} = b_0(-\beta) = \left(1 - \frac{q^2}{1-q^4}\right) q$
and $min\{s(\varphi)\} = s(-2\beta) = q^2$.

The condition $max\{b_0(\varphi)\} < min\{s(\varphi)\} \Leftrightarrow \left(1 - \frac{q^2}{1-q^4}\right) q < q^2 \Leftrightarrow 1 - q - q^2 - q^4 + q^5 < 0$

Numerical analysis implies that the inequality holds for $0.6 < q < 1$.
This is given for $90° \leq \alpha \leq 108°$.

**n = 3** Proof for $\overline{P_{3,3}P_{3,4}}$

$P_{0,0}$, $P_{3,3}$ and $P_{3,4}$ are collinear (see Figure 17) and $b_0(-\beta) = \left(1 - \frac{q^2}{1-q^4}\right) q < q - q^3 < q$

| | |
|---|---|
| $4 \leq n \leq 5$ | $-180° - \beta \geq \varphi_n \geq -180° - \beta - \alpha$ and $r_n \leq \frac{q^3}{1-q^4}$. |
| | Therefore $\overline{P_{n,2^{n-1}-1}P_{n,2^{n-1}}} \subset T_2$ |
| $6 \leq n \leq \frac{360°}{\beta} + 1$ | $-360° + \beta \geq \varphi_n \geq -540° - \beta$ and $r_n \leq e_0(\varphi_n + 360°) = \frac{q^4}{1-q^4} q^{-\frac{\varphi_n}{\beta}}$. |
| | Therefore $\overline{P_{n,2^{n-1}-1}P_{n,2^{n-1}}} \subset S_4 \cup T_1$ |
| $\frac{360°}{\beta} + 1 < n$ | The argument repeats with smaller radii $r_n$. □ |

Proof for $\overline{P_{n,2^{n-1}}P_{n,2^{n-1}+1}}$, $n \geq 2$

With respect to $P_{1,1}$ the points $P_{n,2^{n-1}+1}$ have polar angle of $\varphi_n = -360° - (n-4)\beta$
and distance $r_n = q^n$.

**n = 2** Proof (for $\overline{P_{2,2}P_{2,3}}$)

$90° \leq \alpha \leq 108°$ yields $1 - q^2 \leq q$, so for all points P on segment $\overline{P_{2,2}P_{2,3}}$ the inequality $|P_{0,0}P| \leq q$ holds.
Let F be the foot of the perpendicular from $P_{0,0}$ to the line $P_{2,2}P_{2,3}$. For $90° < \alpha < 120°$ F is an interior point on $\overline{P_{2,2}P_{2,3}}$ and $|P_{0,0}F|$ is the minimum distance from $P_{0,0}$ to $\overline{P_{2,2}P_{2,3}}$.

$|P_{0,0}F| = q \cdot \sin(\alpha - \beta) = \frac{\sqrt{4q^2-1}}{2q^2}(1 - q^2)$

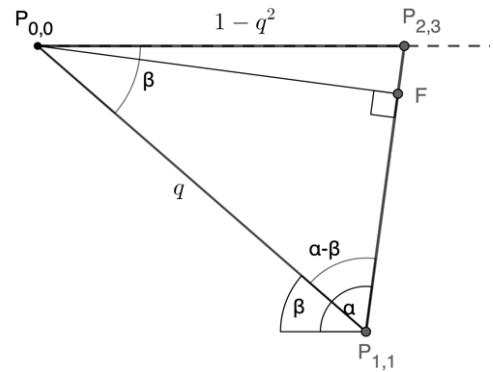

Figure 19: The segment $\overline{P_{2,2}P_{2,3}}$ (note: $P_{2,2} = P_{1,1}$)

$max\{b_0(\varphi): -\beta \leq \varphi \leq 0°\} = b_0(0°) = \left(1 - \frac{q^2}{1-q^4}\right)$.

The condition $|P_{0,0}F| > max\{b_0(\varphi)\}$ is equivalent to $\frac{\sqrt{4q^2-1}}{2q^2}(1-q^2) > \left(1 - \frac{q^2}{1-q^4}\right)$ and can be transformed to $\sqrt{4q^2-1}(1-q^2-q^4+q^6) - 2(q^2-q^4-q^6) > 0$.

Numerical analysis implies that this holds for $q > 0.599$ corresponding to $\alpha \leq 113°$.

Therefore, $\overline{P_{2,2}P_{2,3}} \subset S_2$ for $\frac{1}{2}\sqrt{2} \geq q \geq \frac{\sqrt{5}-1}{2}$.

| | | |
|---|---|---|
| $3 \leq n \leq 5 + \frac{180°}{\beta}$ | $-360° + \beta \geq \varphi_n \geq -540° - \beta$ and $r_n \leq e_0(\varphi_n + 360°) = \frac{q^4}{1-q^4} q^{-\frac{\varphi_n}{\beta}}$. | |
| | Therefore, $\overline{P_{n,2^{n-1}}P_{n,2^{n-1}+1}} \subset S_4 \cup T_1$ | |
| $5 + \frac{180°}{\beta} \leq n \leq 5 + \frac{180°+\alpha}{\beta}$ | $-540° - \beta \geq \varphi_n \geq -540° - \beta - \alpha$ and $r_n \leq \frac{q^3}{1-q^4}$. | |
| | Therefore, $\overline{P_{n,2^{n-1}}P_{n,2^{n-1}+1}} \subset T_2$ | |
| $5 + \frac{180°+\alpha}{\beta} < n$ | The argument repeats with smaller radii $r_n$. □ | |

Lemma 10

For all paperfolding polygons $Q_n, n \in \mathbb{N}, n \geq 2$ the segments $\overline{P_{n,k}P_{n,k+1}}$, $k \epsilon \{1,2, \ldots, 2^n - 2\}$ are subsets of $A_{PF}$.

Remark

The segments $\overline{P_{n,0}P_{n,1}}$ (first segment) and $\overline{P_{n,2^n-1}P_{n,2^n}}$ (last segment) are excluded.

Proof (by mathematical induction)

Base step: $n = 2$ According to Lemma 9 $\overline{P_{2,1}P_{2,2}}$ and $\overline{P_{2,2}P_{2,3}}$ are subsets of $A_{PF}$.

Induction step: Let $\overline{P_{n,k}P_{n,k+1}}$, $k \epsilon \{1,2, \ldots, 2^n - 2\}$ be subset of $A_{PF}$.

Then $\pi_0(\overline{P_{n,k}P_{n,k+1}})$ is subset of $\pi_0(A_{PF})$. $\pi_0(A_{PF}) \subset A_{PF}$ (Proposition 1 a)),

therefore $\pi_0(\overline{P_{n,k}P_{n,k+1}}) = \overline{P_{n+1,k}P_{n+1,k+1}}$, $k \epsilon \{1,2, \ldots, 2^n - 2\}$ is subset of $A_{PF}$.

$\pi_1(\overline{P_{n,k}P_{n,k+1}})$ is subset of $\pi_1(A_{PF})$. $\pi_1(A_{PF}) \subset A_{PF}$ (Proposition 1 b),

therefore $\pi_1(\overline{P_{n,k}P_{n,k+1}}) = \overline{P_{n+1,2^{n+1}-(k+1)}P_{n+1,2^{n+1}-k}}$, $k \epsilon \{1,2, \ldots, 2^n - 2\}$ is subset of $A_{PF}$.

For $k = 2^n - 1$, the segments $\overline{P_{n+1,2^n-1}P_{n+1,2^n}}$ and $\overline{P_{n+1,2^n}P_{n+1,2^n+1}}$ are subsets of $A_{PF}$, as proven in Lemma 9.

Lemma 8 and Lemma 10 yield the proof of Proposition 2.

Now everything is ready for the main sentences.

## 7. Avoiding intersections of the paperfolding polygons $Q_n$

Proposition 3

For unfolding angles $98.195° \leq \alpha \leq 108°$ $\pi_0(A_{PF}) \cap \pi_1(A_{PF}) = \{P_{1,1}\}$.

Proof

The main condition is $c_1(\varphi) < d_2(\varphi)$

$$\Leftrightarrow \frac{1}{1-q^4} q^{-\frac{\varphi+\alpha}{\beta}} < \left(1 - \frac{q^2}{1-q^4}\right) q^{-\frac{\varphi+2\alpha}{\beta}}$$

$$\Leftrightarrow q^{-\frac{\varphi+\alpha}{\beta}} < (1-q^2-q^4) q^{-\frac{\varphi+2\alpha}{\beta}}$$

$$\Leftrightarrow 0 < 1 - q^2 - q^4 - q^{\frac{\alpha}{\beta}}$$

The right-hand side is in [0,1] monotonously decreasing with a zero between $q$ = 0.6615289 and $q$ = 0.6615339. Therefore, $c_1(\varphi) < d_2(\varphi)$ holds for $q \leq 0.6615289$, which corresponds to $\alpha \geq 98.195°$.

To prove that $\pi_0(A_{PF}) \cap \pi_1(A_{PF}) = \{P_{1,1}\}$ for $98.195° \leq \alpha \leq 108°$ in detail, the plane is partitioned into three pairs of half-planes:
- by line $P_{0,0}P_{1,1}$ into $H_0$ containing $P_{2,1}$ and $H_1$ containing $P_{2,3}$
- by line $P_{2,3}P_{1,1}$ into $H_2$ containing $P_{0,1}$ and $H_3$ containing $P_{0,0}$
- by line $P_{0,0}P_{0,1}$ into $H_4$ containing $P_{1,1}$ and $H_5$ not containing $P_{1,1}$

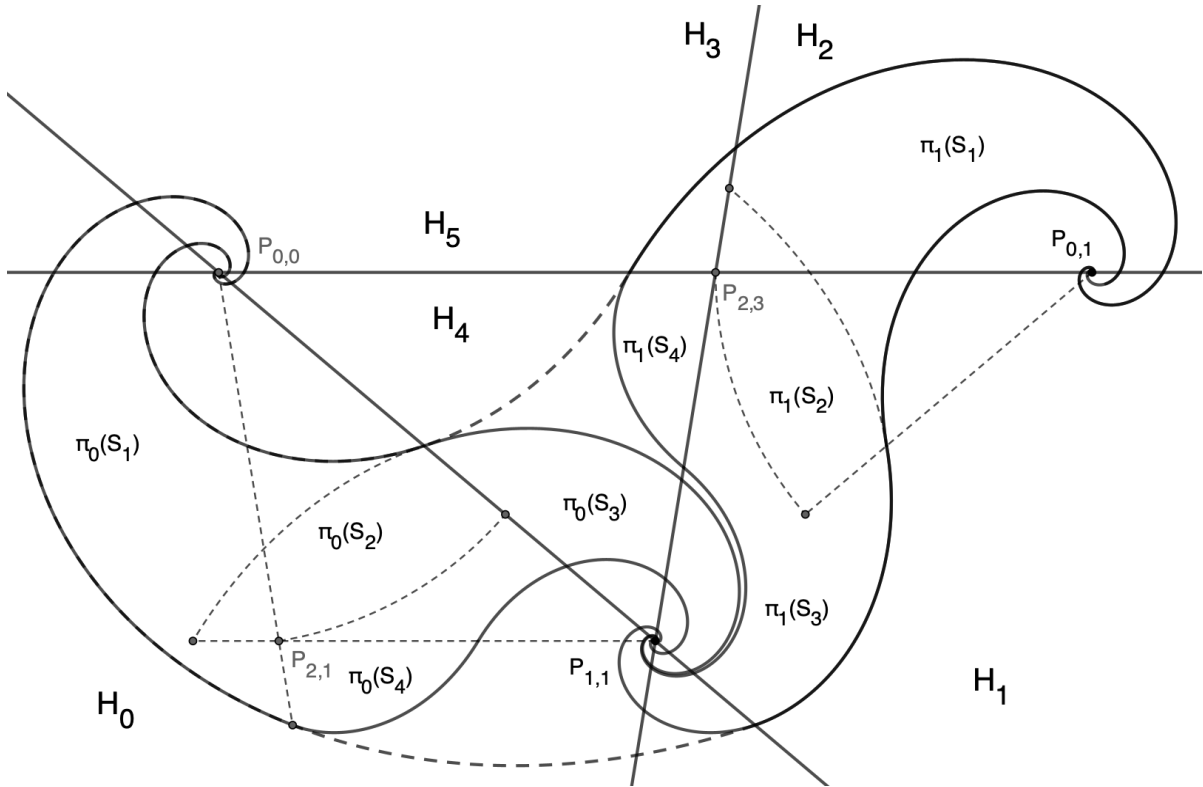

Figure 20: $\alpha = 99°$ $\pi_0(A_{PF})$ and $\pi_1(A_{PF})$ are clearly separated

The poof is organized in the following table:

| ∩ | $\pi_0(S_1)$ | $\pi_0(S_2)$ | $\pi_0(S_3)$ | $\pi_0(S_4)$ |
|---|---|---|---|---|
| $\pi_1(S_1)$ | 1. | ∅ | 2. | ∅ |
| $\pi_1(S_2)$ | 3. | ∅ | 4. | ∅ |
| $\pi_1(S_3)$ | 5. | 6. | 7. | 8. |
| $\pi_1(S_4)$ | 9. | ∅ | 10. | ∅ |

$\pi_0(S_2)$ and $\pi_0(S_4)$ are subsets of $H_0$. $\pi_1(S_1), \pi_1(S_2)$ and $\pi_1(S_4)$ are subsets of $H_1$, because the distance between the point $P_{0,1}$ and the line $P_{0,0}P_{1,1}$ is greater than that between the point $P_{0,1}$ and any point of these three sets. This yields the ∅ in the table.

1. $\pi_0(S_1) \cap \pi_1(S_1)$

   Since $\pi_1(S_1) \subset H_1$, only $\pi_0(S_1) \cap H_1$ needs to be considered in the calculation.
   $$P \in \pi_0(S_1) \cap H_1 \Rightarrow |P_{0,0}P| \leq a_1(-180° - \beta) = \frac{1}{1-q^4} q^{-\frac{-180°-\beta}{\beta}}$$
   Applying the triangle inequality: $|P_{0,1}P| \geq |P_{0,0}P_{0,1}| - |P_{0,0}P| \geq 1 - \frac{1}{1-q^4} q^{-\frac{-180°-\beta}{\beta}}$.
   $$Q \in \pi_1(S_1) \Rightarrow |P_{0,1}Q| \leq a_2(-180°) = \frac{1}{1-q^4} q^{-\frac{-180°+\alpha}{\beta}}$$
   $\frac{1}{1-q^4} q^{-\frac{-180°+\alpha}{\beta}} < 1 - \frac{1}{1-q^4} q^{-\frac{-180°-\beta}{\beta}}$ holds for $q \leq \frac{1}{2}\sqrt{2}$.
   For any $P \in \pi_0(S_1) \cap H_1$ and $Q \in \pi_1(S_1)$, the inequality $|P_{0,1}Q| < |P_{0,1}P|$ holds.
   This implies $(\pi_0(S_1) \cap H_1) \cap \pi_1(S_1) = \emptyset$ and therefore $\pi_0(S_1) \cap \pi_1(S_1) = \emptyset$.

2. $\pi_0(S_3) \cap \pi_1(S_1)$

   Since $\pi_1(S_1) \subset H_1$ only $\pi_0(S_3) \cap H_1$ needs to be considered in the calculation.
   $\pi_1(S_1)$ is split by line $P_{0,0}P_{0,1}$ into $\pi_1(S_1) \cap H_5$ and $\pi_1(S_1) \cap H_4$.
   $P \in \pi_0(S_3) \cap H_1 \Rightarrow |P_{1,1}P| \leq c_1(-180° - \beta) = \frac{1}{1-q^4} q^3$ Lemma 6 a) implies
   that $\pi_0(S_3) \cap H_1 \subset H_4$. So $(\pi_0(S_3) \cap H_1) \cap (\pi_1(S_1) \cap H_5) = \emptyset$.
   Applying the triangle inequality: $|P_{0,1}P| \geq |P_{1,1}P_{0,1}| - |P_{1,1}P| \geq q - \frac{1}{1-q^4} q^3$
   $$Q \in \pi_1(S_1) \cap H_4 \Rightarrow |P_{0,1}Q| \leq a_2(-360°) = \frac{1}{1-q^4} q^{-\frac{-360°+\alpha}{\beta}}$$
   $\frac{1}{1-q^4} q^{-\frac{-360°+\alpha}{\beta}} < q - \frac{1}{1-q^4} q^3$ holds for $q \leq \frac{1}{2}\sqrt{2}$.
   For any $P \in \pi_0(S_3) \cap H_1$ and $Q \in \pi_1(S_1) \cap H_4$, the inequality $|P_{0,1}Q| < |P_{0,1}P|$ holds.
   This implies $(\pi_0(S_3) \cap H_1) \cap (\pi_1(S_1) \cap H_4) = \emptyset$ and therefore $\pi_0(S_3) \cap \pi_1(S_1) = \emptyset$.

3. $\pi_0(S_1) \cap \pi_1(S_2)$

   $\pi_0(S_1) \subset H_3$ and $\pi_1(S_2) \subset H_2$ and as $H_2$ and $H_3$ are disjoint halfplanes with respect to line $P_{2,3}P_{1,1}$ this implies $\pi_0(S_1) \cap \pi_1(S_2) = \emptyset$

4. $\pi_0(S_3) \cap \pi_1(S_2)$

   Since $\pi_1(S_2) \subset H_2$, only $\pi_0(S_3) \cap H_2$ needs to be considered in the calculation.
   $$P \in \pi_0(S_3) \cap H_2 \Rightarrow |P_{1,1}P| \leq c_1(-180° - \alpha) = \frac{1}{1-q^4} q^{-\frac{-180°}{\beta}}$$
   Applying the triangle inequality: $|P_{0,1}P| \geq |P_{1,1}P_{0,1}| - |P_{1,1}P| \geq q - \frac{1}{1-q^4} q^{\frac{\alpha+2\beta}{\beta}}$
   $Q \in \pi_1(S_2) \Rightarrow |P_{0,1}Q| \leq q^2$ and $q^2 < q - \frac{1}{1-q^4} q^{\frac{\alpha+2\beta}{\beta}}$ holds for $q < 0.67$.
   For any $P \in \pi_0(S_3) \cap H_2$ and $Q \in \pi_1(S_2)$, the inequality $|P_{0,1}Q| < |P_{0,1}P|$ holds.
   This implies $(\pi_0(S_3) \cap H_2) \cap \pi_1(S_2) = \emptyset$ and therefore $\pi_0(S_3) \cap \pi_1(S_2) = \emptyset$.

5. $\pi_0(S_1) \cap \pi_1(S_3)$

   The line $P_{0,0}P_{1,1}$ divides both sets. Thus, there are two cases.
   Case 1: $(\pi_0(S_1) \cap H_0) \cap (\pi_1(S_3) \cap H_0)$
   The line $P_{0,0}P_{2,1}$ partitions the plane into two half-planes. $\pi_1(S_3) \cap H_0$ is entirely in one half-plane, say $H'$. Hence only $\pi_0(S_1) \cap H_0 \cap H'$ needs to be considered in the calculation.
   Let P be a point of $\pi_0(S_1) \cap H_0 \cap H'$ with a polar angel $\varphi$ to the center $P_{0,0}$,
   $-\beta - k \cdot 360° \geq \varphi \geq -2\beta - k \cdot 360°, k = 1,2,3,...$.
   Then $|P_{0,0}P| \leq a_1(-\beta - 360°) = \frac{1}{1-q^4} q^{-\frac{-5\beta-2\alpha}{\beta}}$.

Applying the triangle inequality: $|P_{1,1}P| \geq |P_{0,0}P_{1,1}| - |P_{0,0}P| \leq q - \frac{1}{1-q^4}q^{-\frac{-5\beta-2\alpha}{\beta}}$

$Q \in \pi_1(S_3) \cap H_0 \Rightarrow |P_{1,1}Q| \leq c_2(-\beta - 360°) = \frac{1}{1-q^4}q^{-\frac{-\beta-360°+2\alpha}{\beta}}$

$\frac{1}{1-q^4}q^{-\frac{-\beta-360°+2\alpha}{\beta}} < q - \frac{1}{1-q^4}q^{-\frac{-5\beta-2\alpha}{\beta}}$ holds for $q \leq \frac{1}{2}\sqrt{2}$.

For any $P \in \pi_0(S_1) \cap H_0 \cap H'$ and $Q \in \pi_1(S_3) \cap H_0$, the inequality $|P_{1,1}Q| < |P_{1,1}P|$ holds.
This implies $(\pi_0(S_1) \cap H_0 \cap H') \cap (\pi_1(S_3) \cap H_0) = \emptyset$
and therefore $(\pi_0(S_1) \cap H_0) \cap (\pi_1(S_3) \cap H_0) = \emptyset$.

Case 2: $(\pi_0(S_1) \cap H_1) \cap (\pi_1(S_3) \cap H_1)$
Since $\pi_0(S_1) \cap H_1 \subset H_3$ only $\pi_1(S_3) \cap H_1 \cap H_3$ needs to be considered in the calculation.

$P \in \pi_0(S_1) \cap H_1 \Rightarrow |P_{0,0}P| \leq a_1(-180° - \beta) = \frac{1}{1-q^4}q^{-\frac{-180°-\beta}{\beta}}$.

Applying the triangle inequality: $|P_{1,1}P| \geq |P_{0,0}P_{1,1}| - |P_{0,0}P| \geq q - \frac{1}{1-q^4}q^{-\frac{-180°-\beta}{\beta}}$.

$Q \in \pi_1(S_3) \cap H_1 \cap H_3 \Rightarrow |P_{1,1}Q| \leq c_2(-\beta - 540°) = \frac{1}{1-q^4}q^{-\frac{-\beta-540°+2\alpha}{\beta}}$

$\frac{1}{1-q^4}q^{-\frac{-\beta-540°+2\alpha}{\beta}} < q - \frac{1}{1-q^4}q^{-\frac{-180°-\beta}{\beta}}$ holds for $q \leq \frac{1}{2}\sqrt{2}$.

For any $P \in \pi_0(S_1) \cap H_1$ and $Q \in \pi_1(S_3) \cap H_1 \cap H_3$, the inequality $|P_{1,1}Q| < |P_{1,1}P|$ holds.
This implies $(\pi_0(S_1) \cap H_1) \cap (\pi_1(S_3) \cap H_1 \cap H_3) = \emptyset$
and therefore $(\pi_0(S_1) \cap H_1) \cap (\pi_1(S_3) \cap H_1) = \emptyset$.

6. $\pi_0(S_2) \cap \pi_1(S_3)$
   Since $\pi_0(S_2) \subset H_0$ only $\pi_1(S_3) \cap H_0$ needs to be considered in the calculation.
   $P \in \pi_0(S_2) \Rightarrow |P_{0,0}P| \leq q^2$
   Applying the triangle inequality: $|P_{1,1}P| \geq |P_{0,0}P_{1,1}| - |P_{0,0}P| \geq q - q^2$
   $Q \in \pi_1(S_3) \cap H_0 \Rightarrow |P_{1,1}Q| \leq c_2(-\beta - 360°) = \frac{1}{1-q^4}q^{-\frac{-\beta-360°+2\alpha}{\beta}}$
   $\frac{1}{1-q^4}q^{-\frac{-\beta-360°+2\alpha}{\beta}} < q - q^2$ holds for $q < 0.69$.
   For any $P \in \pi_0(S_2)$ and $Q \in \pi_1(S_3) \cap H_0$, the inequality $|P_{1,1}Q| < |P_{1,1}P|$ holds.
   This implies $\pi_0(S_2) \cap (\pi_1(S_3) \cap H_0) = \emptyset$ and therefore $\pi_0(S_2) \cap \pi_1(S_3) = \emptyset$.

7. $\pi_0(S_3) \cap \pi_1(S_3)$
   For polar angles $-180° \geq \varphi > -180° - \alpha$ only $\pi_0(S_3)$ is defined.
   For $-180° - \alpha \geq \varphi$ a ray from $P_{1,1}$ at polar angle $\varphi$ intersects $\pi_0(S_3)$ and $\pi_1(S_3)$ infinitely many times. The situation is illustrated in Figure 21.

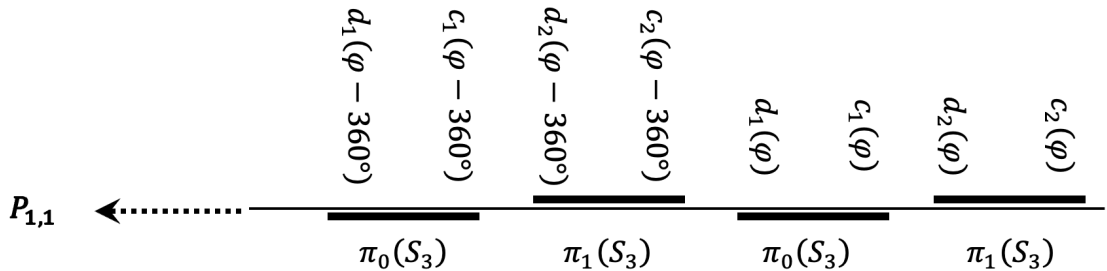

Figure 21: The infinitely intersection of $\pi_0(S_3)$ and $\pi_1(S_3)$ (case 7)

Points between $c_2(\varphi - k \cdot 360°)$ and $d_2(\varphi - k \cdot 360°)$, $k \in \mathbb{N}_0$ belong to $\pi_1(S_3)$,
points between $c_1(\varphi - k \cdot 360°)$ and $d_1(\varphi - k \cdot 360°)$, $k \in \mathbb{N}_0$ belong to $\pi_0(S_3)$.

These intervals are pairwise disjoint if $c_2(\varphi) > d_2(\varphi) > c_1(\varphi) > d_1(\varphi) > c_2(\varphi - 360°)$ for all $-180° - \alpha \geq \varphi$.

$\frac{c_2(\varphi)}{d_2(\varphi)} = \frac{1}{1-q^4-q^2} > 1$ holds for $q \leq \frac{1}{2}\sqrt{2}$

$d_2(\varphi) > c_1(\varphi)$ was calculated in the beginning,
yielding the condition $q \leq 0.6615289$, which corresponds to $\alpha \geq 98.195°$.

$\frac{c_1(\varphi)}{d_1(\varphi)} = \frac{1}{1-q^4-q^2} > 1$ holds for $q \leq \frac{1}{2}\sqrt{2}$

$d_1(\varphi) = \left(1 - \frac{q^2}{1-q^4}\right) q^{-\frac{\varphi+\alpha}{\beta}}$, $c_2(\varphi - 360°) = \frac{1}{1-q^4} q^{-\frac{\varphi - 4\beta}{\beta}}$

$d_1(\varphi) > c_2(\varphi - 360°)$ simplifies to $1 - q^2 - q^4 - q^{4+\frac{\alpha}{\beta}} > 0$, which holds for $q \leq \frac{1}{2}\sqrt{2}$.

Therefore, $\pi_0(S_3) \cap \pi_1(S_3) = \{P_{1,1}\}$.

8. $\pi_0(S_4) \cap \pi_1(S_3)$
   Since $\pi_0(S_4) \subset H_0$ only $\pi_1(S_3) \cap H_0$ needs to be considered in the calculation.
   $P \in \pi_0(S_4) \Rightarrow |P_{2,1}P| \leq e_1(0°) = \frac{1}{1-q^4} q^{-\frac{4\beta}{\beta}}$.
   
   Applying the triangle inequality: $|P_{1,1}P| \geq |P_{2,1}P_{1,1}| - |P_{2,1}P| \geq q^2 - \frac{1}{1-q^4} q^{-\frac{4\beta}{\beta}}$.
   
   $Q \in \pi_1(S_3) \cap H_0 \Rightarrow |P_{1,1}Q| \leq c_2(-\beta - 360°) = \frac{1}{1-q^4} q^{-\frac{-\beta - 360° + 2\alpha}{\beta}}$
   
   $\frac{1}{1-q^4} q^{-\frac{-\beta - 360° + 2\alpha}{\beta}} < q^2 - \frac{1}{1-q^4} q^{-\frac{4\beta}{\beta}}$ holds for $q < 0.68$.
   
   For any $P \in \pi_0(S_4)$ and $Q \in \pi_1(S_3) \cap H_0$, the inequality $|P_{1,1}Q| < |P_{1,1}P|$ holds.
   This implies $\pi_0(S_4) \cap (\pi_1(S_3) \cap H_0) = \emptyset$ and therefore $\pi_0(S_4) \cap \pi_1(S_3) = \emptyset$.

9. $\pi_0(S_1) \cap \pi_1(S_4)$
   Since $\pi_1(S_4) \subset H_1$ only $\pi_0(S_1) \cap H_1$ needs to be considered in the calculation.
   $P \in \pi_0(S_1) \cap H_1 \Rightarrow |P_{0,0}P| \leq a_1(-180° - \beta) = \frac{1}{1-q^4} q^{-\frac{-180° - \beta}{\beta}}$.
   
   Applying the triangle inequality: $|P_{2,3}P| \geq |P_{0,0}P_{2,3}| - |P_{0,0}P| \geq 1 - q^2 - \frac{1}{1-q^4} q^{-\frac{-180° - \beta}{\beta}}$.
   
   $Q \in \pi_1(S_4) \Rightarrow |P_{2,3}Q| \leq e_2(-\alpha) = \frac{1}{1-q^4} q^{-\frac{4\beta}{\beta}}$
   
   $\frac{1}{1-q^4} q^{-\frac{4\beta}{\beta}} < 1 - q^2 - \frac{1}{1-q^4} q^{-\frac{-180° - \beta}{\beta}}$ holds for $q < 0.69$.
   
   For any $P \in \pi_0(S_1) \cap H_1$ and $Q \in \pi_1(S_4)$, the inequality $|P_{2,3}Q| < |P_{2,3}P|$ holds.
   This implies $(\pi_0(S_1) \cap H_1) \cap \pi_1(S_4) = \emptyset$ and therefore $\pi_0(S_1) \cap \pi_1(S_4) = \emptyset$.

10. $\pi_0(S_3) \cap \pi_1(S_4)$
    $c_1(-180° - \alpha)$ defines the point $R_1$ on the spiral $c_1$ and $e_2(-\alpha)$ defines the point $R_2$ on the spiral $e_2$.
    Let $t_1$ be the tangent line to $c_1$ at $R_1$.
    Let $t_2$ be the tangent line to $e_2$ at $R_2$.
    Since $P_{1,1}$, $R_1$, $R_2$ and $P_{2,3}$ are collinear, $t_1$ and $t_2$ are parallel.
    $c_1(-180° - \alpha) < d_2(-180° - \alpha) = |P_{2,3}P_{1,1}| - e_2(-\alpha)$ implies that $t_1$ and $t_2$ are distinct lines.
    $c_1$ and $t_2$ lie in different half-planes determined by $t_1$.
    $e_2$ and $t_1$ lie in different half-planes determined by $t_2$.
    Therefore, $\pi_0(S_3) \cap \pi_1(S_4) = \emptyset$.

    □

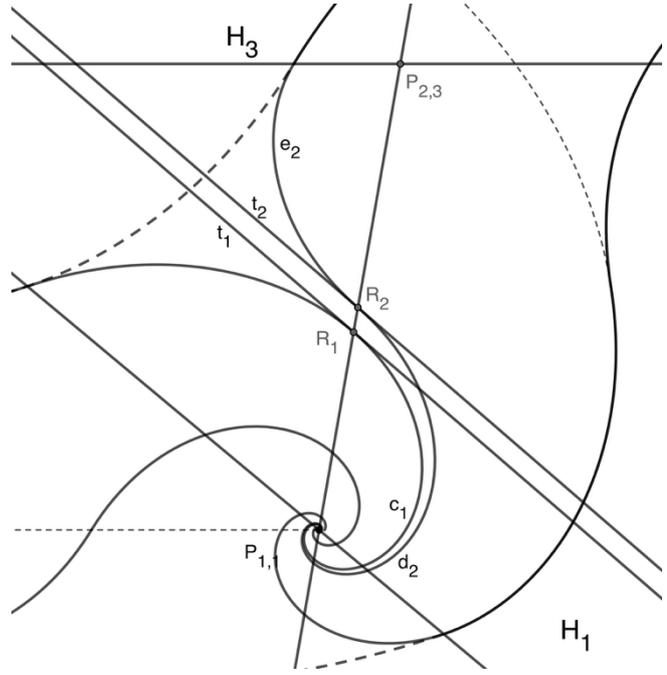

Figure 22: $\pi_0(S_3)$ and $\pi_1(S_4)$ (case 10)

Corollary 2 (for generalized $A_{PF}$)
For unfolding angles $98.195° \leq \alpha \leq 108°$
$A_{PF}\left(\overline{P_{n,2k \cdot 2^m} P_{n,(2k+1) \cdot 2^m}}\right) \cap A_{PF}\left(\overline{P_{n,(2k+2) \cdot 2^m} P_{n,(2k+1) \cdot 2^m}}\right) = \{P_{n,(2k+1) \cdot 2^m}\}$.

Theorem 1
For unfolding angles $98.195° \leq \alpha \leq 108°$ the paperfolding polygons $Q_n$ avoid intersections.

Proof
The paperfolding polygon $Q_n$ is constructed by $Q_n = \pi_0(Q_{n-1}) \cup \pi_1(Q_{n-1})$.
$\pi_0(A_{PF}) \cap \pi_1(A_{PF}) = \{P_{1,1}\}$ implies that all segments of $\pi_0(Q_{n-1})$, which are subsets of $\pi_0(A_{PF})$ do not intersect any segment of $\pi_1(Q_{n-1})$, which is a subset of $\pi_1(A_{PF})$.
This argument excludes the segment $\overline{P_{n,0}P_{n,1}}$ at the beginning and $\overline{P_{n,2^{n-1}-1}P_{n,2^{n-1}}}$ at the end of $\pi_0(Q_{n-1})$, respectively $\overline{P_{n,2^{n-1}}P_{n,2^{n-1}+1}}$ at the beginning and $\overline{P_{n,2^n-1}P_{n,2^n}}$ at the end of $\pi_1(Q_{n-1})$.
Discussing theses four segments:
1. $\overline{P_{n,0}P_{n,1}}$ does not intersect with any segment of $\pi_1(Q_{n-1})$, as for $n \geq 3$ the circle with center $P_{0,0}$ and radius $q^n$ does not intersect $\pi_1(A_{PF})$.

2. An analogous argumentation holds for the segment $\overline{P_{n,2^n-1}P_{n,2^n}}$ at the end of $\pi_1(Q_{n-1})$.

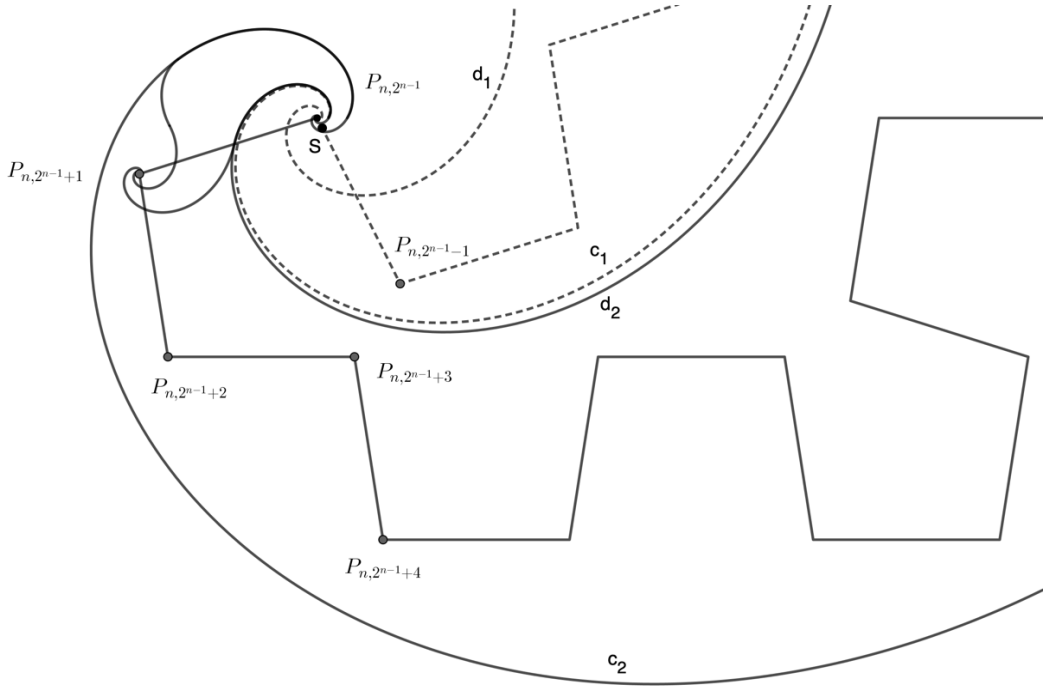

Figure 23: The situation at $P_{n,2^{n-1}}$

3. Assume toward a contradiction that segment $\overline{P_{n,2^{n-1}-1}P_{n,2^{n-1}}}$ intersects in point S a segment $\overline{P_{n,k}P_{n,k+1}}$, which is a segment of $\pi_1(Q_{n-1})$. Therefore, $2^{n-1} \leq k \leq 2^n - 1$.

   Case 1: $S \notin \pi_1(A_{PF})$

   Then S can only be the intersection point with the first segment $\overline{P_{n,2^{n-1}}P_{n,2^{n-1}+1}}$ of $\pi_1(Q_{n-1})$ or with the last segment $\overline{P_{n,2^n-1}P_{n,2^n}}$.

   First segment $\overline{P_{n,2^{n-1}}P_{n,2^{n-1}+1}}$: If $S = \overline{P_{n,2^{n-1}}P_{n,2^{n-1}+1}} \cap \overline{P_{n,2^{n-1}-1}P_{n,2^{n-1}}}$ then $S = P_{n,2^{n-1}}$. This is not an intersection in the discussed sense.

   Last segment $\overline{P_{n,2^n-1}P_{n,2^n}}$: Let Q be a point on $\overline{P_{n,2^n-1}P_{n,2^n}}$. Then $|QP_{n,2^n}| \leq q^n$.

   $S \in \overline{P_{n,2^{n-1}-1}P_{n,2^{n-1}}}$ yields $|SP_{n,2^{n-1}}| \leq q^n$.

   Applying the triangle inequality: $|SP_{n,2^n}| \geq |P_{n,2^n-1}P_{n,2^n}| - |SP_{n,2^{n-1}}| = q - q^n$

   $q - q^n > q^n$ holds for $q \leq \frac{1}{2}\sqrt{2}$ and $n \geq 4$.

   For $S \in \overline{P_{n,2^{n-1}-1}P_{n,2^{n-1}}}$ and any $Q \in \overline{P_{n,2^n-1}P_{n,2^n}}$, the inequality $|QP_{n,2^n}| < |SP_{n,2^n}|$ holds. This implies $S \notin \overline{P_{n,2^n-1}P_{n,2^n}}$ and therefore S is not the intersection point of $\overline{P_{n,2^n-1}P_{n,2^n}}$ and $\overline{P_{n,2^{n-1}-1}P_{n,2^{n-1}}}$.

   Case 2: $S \in \pi_1(A_{PF})$

   $S \in \overline{P_{n,2^{n-1}-1}P_{n,2^{n-1}}}$ and $\overline{P_{n,2^{n-1}-1}P_{n,2^{n-1}}}$ has a polar angle of $-180° - (n-2)\beta$, with respect to the point $P_{n,2^{n-1}}$. Therefore

   $$|SP_{n,2^{n-1}}| \leq c_2(-540° - (n-2)\beta) = \frac{1}{1-q^4}q^{-\frac{-540°-(n-2)\beta+2\alpha}{\beta}} = \frac{1}{1-q^4}q^n q^4 q^{\frac{\alpha}{\beta}}$$

   Hence $\frac{1}{1-q^4}q^n q^4 q^{\frac{\alpha}{\beta}} < q^n$ the point S is an element of $A_{PF}\left(\overline{P_{n,2^{n-1}}P_{n,2^{n-1}+1}}\right)$ which implies that S belongs to all $A_{PF}\left(\overline{P_{n,2^{n-1}}P_{n,2^{n-1}+2m}}\right)$, $0 \leq m \leq n-2$ (Corollary 1 of

Proposition 1).

Now the main idea is: Because $S \in A_{PF}(\overline{P_{n,2^{n-1}}P_{n,2^{n-1}+2^m}})$, $0 \leq m \leq n-2$ and
$A_{PF}(\overline{P_{n,2^{n-1}}P_{n,2^{n-1}+2^m}}) \cap A_{PF}(\overline{P_{n,2^{n-1}+2^{m+1}}P_{n,2^{n-1}+2^m}}) = \{P_{n,2^{n-1}+2^m}\}$ (Corollary 2)
and $S \neq P_{n,2^{n-1}+2^m}$ (not an intersection in the discussed sense),
it follows that S cannot be an element of $A_{PF}(\overline{P_{n,2^{n-1}+2^{m+1}}P_{n,2^{n-1}+2^m}})$.
Therefore, S cannot be the intersection point between the segment $\overline{P_{n,2^{n-1}-1}P_{n,2^{n-1}}}$
and any segment $\overline{P_{n,2^{n-1}+k}P_{n,2^{n-1}+k+1}}$, $2^m + 1 \leq k \leq 2^{m+1} - 2$, i.e. any segment that
is a subset of $A_{PF}(\overline{P_{n,2^{n-1}+2^{m+1}}P_{n,2^{n-1}+2^m}})$. The special segments
$\overline{P_{n,2^{n-1}+2^m}P_{n,2^{n-1}+2^m+1}}$ (first) and $\overline{P_{n,2^{n-1}+2^{m+1}-1}P_{n,2^{n-1}+2^{m+1}}}$ (last) have to be treated
separately.

The proof proceeds by case analysis on the parameter $m$:

$m = 0$:

Since $S \in A_{PF}(\overline{P_{n,2^{n-1}}P_{n,2^{n-1}+1}})$ it follows $S \notin A_{PF}(\overline{P_{n,2^{n-1}+1}P_{n,2^{n-1}+2}})$

The special case of $A_{PF}(\overline{P_{n,2^{n-1}+1}P_{n,2^{n-1}+2}})$ concerns the segment $\overline{P_{n,2^{n-1}+1}P_{n,2^{n-1}+2}}$.
All points of this line segment have polar angles $\varphi$ with respect to $P_{n,2^{n-1}}$ satisfying
$-180° - (n-2)\beta - \alpha + \beta \geq \varphi \geq -180° - (n-2)\beta - \alpha$. S lies on segment
$\overline{P_{n,2^{n-1}-1}P_{n,2^{n-1}}}$ which has a polar angle of $\varphi = -180° - (n-2)\beta$ with respect to
$P_{n,2^{n-1}}$. Since these angle ranges are disjoint, S cannot be the intersection point
between both segments.

$m = 1$:

Since $S \in A_{PF}(\overline{P_{n,2^{n-1}}P_{n,2^{n-1}+2}})$ it follows $S \notin A_{PF}(\overline{P_{n,2^{n-1}+2}P_{n,2^{n-1}+4}})$

Only the special segments have to be considered.

First segment $\overline{P_{n,2^{n-1}+2}P_{n,2^{n-1}+3}}$: All points of this line segment have polar angles $\varphi$
with respect to $P_{n,2^{n-1}}$ satisfying
$-180° - (n-2)\beta - \alpha + 2\beta \geq \varphi \geq -180° - (n-2)\beta - \alpha + \beta$. S lies on segment
$\overline{P_{n,2^{n-1}-1}P_{n,2^{n-1}}}$ which has a polar angle of $\varphi = -180° - (n-2)\beta$ with respect to
$P_{n,2^{n-1}}$. Since these angle ranges are disjoint, S cannot be the intersection point
between both segments.

Last segment $\overline{P_{n,2^{n-1}+3}P_{n,2^{n-1}+4}}$: All points of this line segment have a polar angle
$\varphi = -180° - (n-2)\beta - \alpha + 2\beta$ with respect to $P_{n,2^{n-1}}$. S lies on segment
$\overline{P_{n,2^{n-1}-1}P_{n,2^{n-1}}}$ which has a polar angle of $\varphi = -180° - (n-2)\beta$ with respect to
$P_{n,2^{n-1}}$. Since these angles are distinct, S cannot be the intersection point between
both segments.

General case : $2 \leq m \leq n-2$

Since $S \in A_{PF}(\overline{P_{n,2^{n-1}}P_{n,2^{n-1}+2^m}})$ it follows that $S \notin A_{PF}(\overline{P_{n,2^{n-1}+2^m}P_{n,2^{n-1}+2^{m+1}}})$.
Thus, S is not the intersection point between $\overline{P_{n,2^{n-1}-1}P_{n,2^{n-1}}}$ and any segment
$\overline{P_{n,k}P_{n,k+1}}$, $2^{n-1} + 2^m + 1 \leq k \leq 2^{n-1} + 2^{m+1} - 2$.

First segment $\overline{P_{n,2^{n-1}+2^m}P_{n,2^{n-1}+2^m+1}}$:

Let Q be a point on $\overline{P_{n,2^{n-1}+2^m}P_{n,2^{n-1}+2^m+1}}$. Then $|QP_{n,2^{n-1}+2^m}| \leq q^n$.
Applying the triangle inequality to S yields

$$\left|SP_{n,2^{n-1}+2^m}\right| \geq \left|P_{n,2^{n-1}}P_{n,2^{n-1}+2^m}\right| - \left|SP_{n,2^{n-1}}\right| = q^{n-m} - \frac{1}{1-q^4}q^n q^4 q^{\frac{\alpha}{\beta}}$$

$q^{n-m} - \frac{1}{1-q^4}q^n q^4 q^{\frac{\alpha}{\beta}} > q^n$ holds for $q \leq \frac{1}{2}\sqrt{2}$ and $m \geq 2$.

For $S \in \overline{P_{n,2^{n-1}-1}P_{n,2^{n-1}}}$ and any $Q \in \overline{P_{n,2^{n-1}+2^m}P_{n,2^{n-1}+2^m+1}}$, the inequality $\left|QP_{n,2^{n-1}+2^m}\right| < \left|SP_{n,2^{n-1}+2^m}\right|$ holds. This implies $S \notin \overline{P_{n,2^{n-1}+2^m}P_{n,2^{n-1}+2^m+1}}$ and therefore S is not the intersection point of $\overline{P_{n,2^n-1}P_{n,2^n}}$ and $\overline{P_{n,2^{n-1}+2^m}P_{n,2^{n-1}+2^m+1}}$.

Last segment $\overline{P_{n,2^{n-1}+2^{m+1}-1}P_{n,2^{n-1}+2^{m+1}}}$:

Let Q be a point on $\overline{P_{n,2^{n-1}+2^{m+1}-1}P_{n,2^{n-1}+2^{m+1}}}$. Then $\left|QP_{n,2^{n-1}+2^{m+1}}\right| \leq q^n$.

For S the triangle inequality yields

$$\left|SP_{n,2^{n-1}+2^{m+1}}\right| \geq \left|P_{n,2^{n-1}}P_{n,2^{n-1}+2^{m+1}}\right| - \left|SP_{n,2^{n-1}}\right| = q^{n-(m+1)} - \frac{1}{1-q^4}q^n q^4 q^{\frac{\alpha}{\beta}}$$

$q^{n-(m+1)} - \frac{1}{1-q^4}q^n q^4 q^{\frac{\alpha}{\beta}} > q^n$ holds for $q \leq \frac{1}{2}\sqrt{2}$ and $m \geq 2$.

For $S \in \overline{P_{n,2^{n-1}-1}P_{n,2^{n-1}}}$ and any $Q \in \overline{P_{n,2^{n-1}+2^{m+1}-1}P_{n,2^{n-1}+2^{m+1}}}$, the inequality $\left|QP_{n,2^{n-1}+2^{m+1}}\right| < \left|SP_{n,2^{n-1}+2^{m+1}}\right|$ holds. This implies $S \notin \overline{P_{n,2^{n-1}+2^{m+1}-1}P_{n,2^{n-1}+2^{m+1}}}$ and therefore S is not the intersection point of $\overline{P_{n,2^n-1}P_{n,2^n}}$ and $\overline{P_{n,2^{n-1}+2^{m+1}-1}P_{n,2^{n-1}+2^{m+1}}}$.

Therefore, segment $\overline{P_{n,2^{n-1}}P_{n,2^{n-1}-1}}$ does not intersect any segment of $\pi_1(Q_{n-1})$.

4. An analogous argumentation like 3. holds for $\overline{P_{n,2^{n-1}}P_{n,2^{n-1}+1}}$ and its intersection with segments of $\pi_0(Q_{n-1})$.

## Applying the duck to the gap between $P_{4,7}$ and $P_{4,11}$

Each segment of a paperfolding polygon develops by the iterative process into a smaller copy of a paperfolding polygon. Consequently, the local duck-shaped hull of any segment serves as a boundary for all following segments between its start and endpoint.

This principle is applied to the loop formed by the points $P_{4,7}$, $P_{4,8}$, …, $P_{4,11}$ and the gap between $P_{4,7}$ and $P_{4,11}$. If the iterative process of generating the paperfolding polygons avoids intersection in this region, the assumption is, that intersections are avoided for all paperfolding polygons[2].

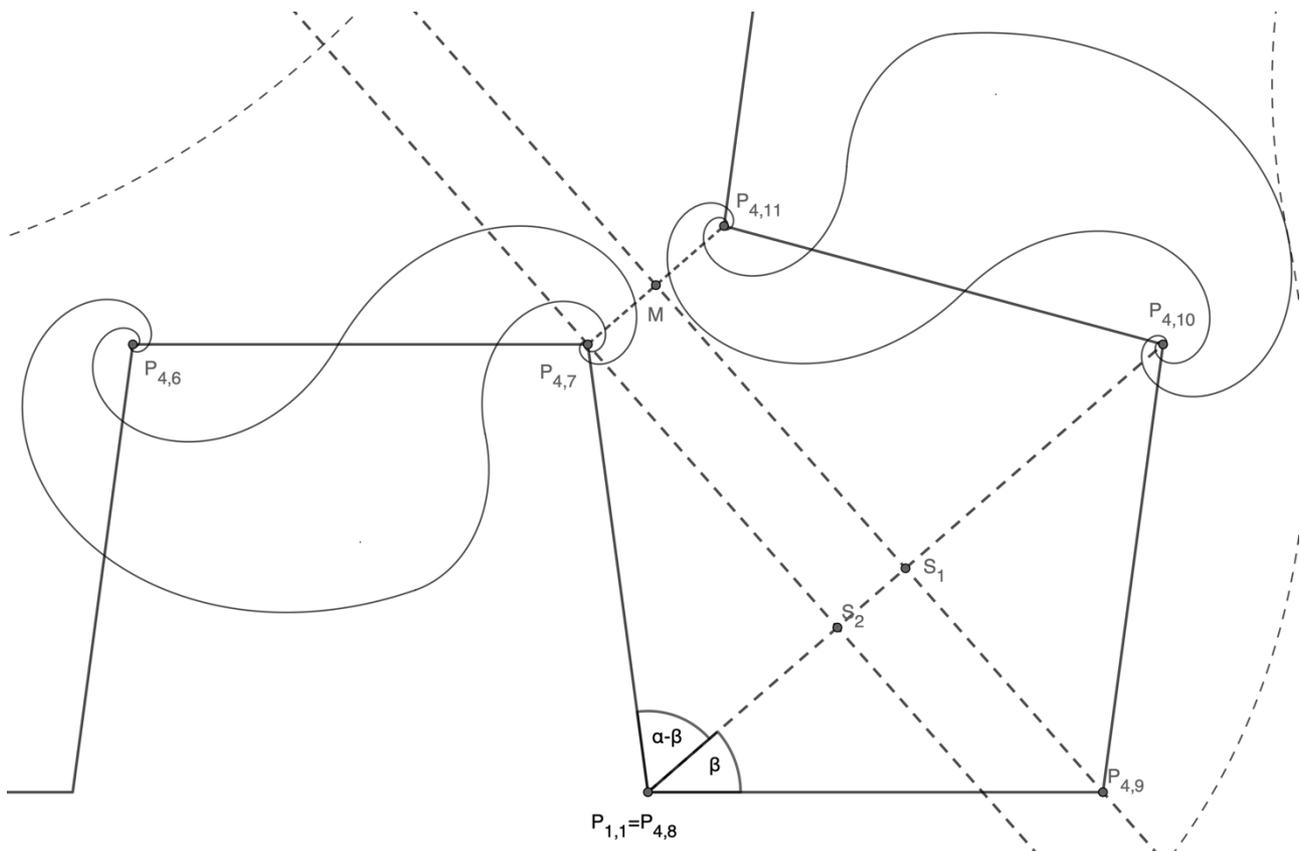

Figure 24: $A_{PF}(\overrightarrow{P_{4,6}P_{4,7}})$ and $A_{PF}(\overrightarrow{P_{4,10}P_{4,11}})$

<u>Lemma 11</u>
If the unfolding angle $\alpha \geq 96.241°$ then $A_{PF}(\overrightarrow{P_{4,6}P_{4,7}})$ and $A_{PF}(\overrightarrow{P_{4,10}P_{4,11}})$ will not intersect.

Proof
The length of all segments in $Q_4$ is $q^4$. For points $S_1$ and $S_2$ see Figure 24.
$|P_{4,7}P_{4,11}| = 2\left(|P_{4,8}S_1| - |P_{4,8}S_2|\right)$
$= 2q^4(cos\beta - cos(\alpha - \beta)) = 2q^4(cos\beta - cos(180° - 3\beta)) = 2q^4(cos\beta + cos3\beta)$

---

[2] D. Knuth writes [1]: „The most troublesome self-crossing points appear to lie near $7 \cdot 2^n$ and $11 \cdot 2^n$"

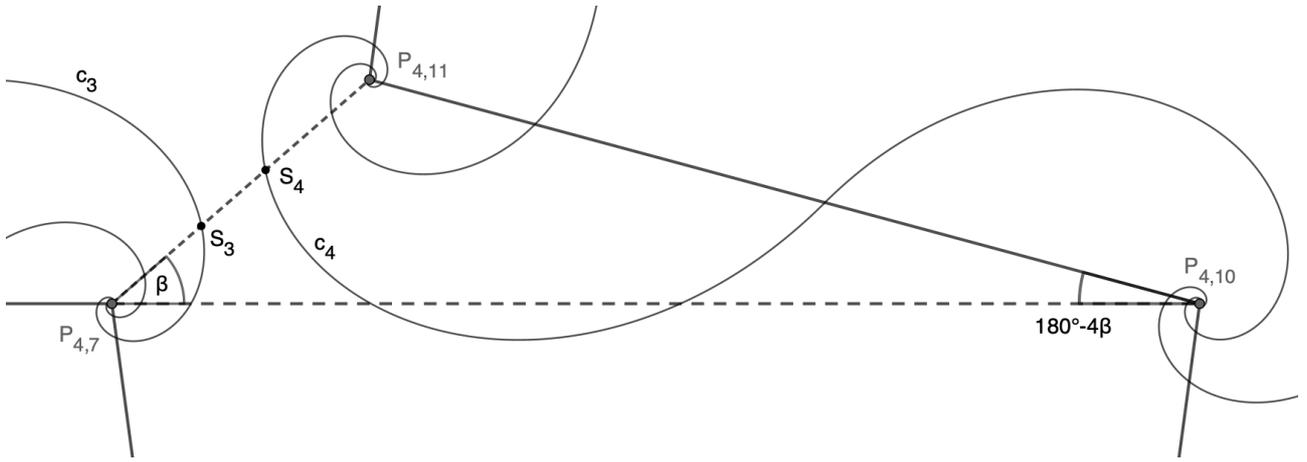

Figure 25: Enlargement of Fig.24

For points $S_3$ and $S_4$ see Figure 25.

An intersection of both ducks is avoided, if $|P_{4,7}S_3| + |P_{4,11}S_4| < |P_{4,7}P_{4,11}|$

Segment $\overline{P_{4,6}P_{4,7}}$ is, compared to $\overline{P_{0,0}P_{0,1}}$, dilated by $q^4$ and not rotated.

So, the outer spiral at $P_{4,7}$ is given by $c_3(\varphi) = q^4 \frac{1}{1-q^4} q^{-\frac{\varphi+\alpha}{\beta}}, \varphi \leq -180°$, center $P_{4,7}$

Therefore, $|P_{4,7}S_3| = c_3(-360° + \beta) = c_3(-2\alpha - 4\beta + \beta) = q^4 \frac{1}{1-q^4} q^{-\frac{-\alpha-3\beta}{\beta}}$.

Segment $\overline{P_{4,10}P_{4,11}}$ is, compared to $\overline{P_{0,0}P_{0,1}}$, dilated by $q^4$ and rotated by $4\beta$.

So, the outer spiral around $P_{4,11}$ is given by

$$c_4(\varphi) = q^4 \frac{1}{1-q^4} q^{-\frac{\varphi+\alpha-4\beta}{\beta}}, \varphi \leq -180° + 4\beta, \quad \text{center } P_{4,11}$$

Therefore, $|P_{4,11}S_4| = c_4(-180° + \beta) = c_4(-\alpha - 2\beta + \beta) = q^4 \frac{1}{1-q^4} q^5$

This yields

$|P_{4,7}S_3| + |P_{4,11}S_4| < |P_{4,7}P_{4,11}|$

$\Leftrightarrow \frac{1}{1-q^4} q^{-\frac{-\alpha-3\beta}{\beta}} + \frac{1}{1-q^4} q^5 < 2(\cos\beta + \cos 3\beta)$

$\Leftrightarrow \frac{q^3}{1-q^4}\left(q^{\frac{\alpha}{\beta}} + q^2\right) < 2(\cos\beta + \cos 3\beta)$

The inequality is solved numerically (3 digits):

| alpha(deg) | beta(deg) | q | left side | right side | < |
|---|---|---|---|---|---|
| 96.235 | 41.883 | 0.6715777 | 0.3238188 | 0.3234373 | FALSE |
| 96.236 | 41.882 | 0.6715724 | 0.3238001 | 0.3234915 | FALSE |
| 96.237 | 41.882 | 0.6715672 | 0.3237814 | 0.3235457 | FALSE |
| 96.238 | 41.881 | 0.6715619 | 0.3237627 | 0.3235999 | FALSE |
| 96.239 | 41.881 | 0.6715567 | 0.3237439 | 0.3236541 | FALSE |
| 96.240 | 41.880 | 0.6715514 | 0.3237252 | 0.3237083 | FALSE |
| 96.241 | 41.880 | 0.6715462 | 0.3237065 | 0.3237625 | TRUE |
| 96.242 | 41.879 | 0.6715409 | 0.3236878 | 0.3238167 | TRUE |
| 96.243 | 41.879 | 0.6715357 | 0.3236691 | 0.3238709 | TRUE |
| 96.244 | 41.878 | 0.6715304 | 0.3236504 | 0.3239251 | TRUE |
| 96.245 | 41.878 | 0.6715252 | 0.3236317 | 0.3239793 | TRUE |

Conjecture

If intersections of $A_{PF}\left(\overrightarrow{P_{4,6}P_{4,7}}\right)$ and $A_{PF}\left(\overrightarrow{P_{4,10}P_{4,11}}\right)$ are avoided in the gap between $P_{4,7}$ and $P_{4,11}$ then intersections are avoided for all paperfolding polygons. By Lemma 11, this is satisfied when the unfolding angle $\alpha \geq 96{,}241°$.

**Summary**

For the unfolding angle $\alpha \leq 90°$ of the paperfolding polygons we know:

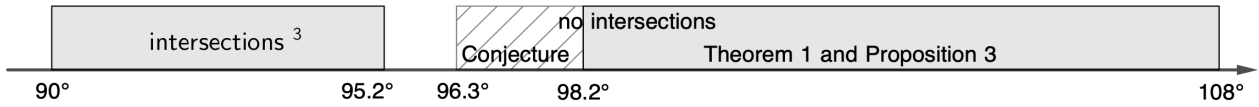

Figure 26: Diagram of the result of this article

| | |
|---|---|
| $90° < \alpha \leq 95.126°$ | Intersections (from $10^{th}$ iteration) [3] |
| $95.126 < \alpha < 96.241$ | Interval of uncertainty |
| $96.241 \leq \alpha < 98.195$ | Conjecture, that there are no intersections |
| $98.195 \leq \alpha \leq 108°$ | No intersections, proved in theorem 1 |

In addition

| | |
|---|---|
| $99.3438 < \alpha \leq 180°$ | No intersections [7] |